\documentclass{amsart}

\author{Euan Aitken}
\address{Radboud University, Nijmegen}
\email{euan.aitken@ru.nl}

\usepackage{array}
\usepackage{amsmath,amssymb,amsthm}
\usepackage[T1]{fontenc}
\usepackage{tgpagella}
\usepackage{graphicx} 
\usepackage{tikz,tikz-cd}
\usepackage{quiver}
\usepackage{hyperref}
\usepackage{enumitem}  
\usepackage{cancel}
\usepackage{adjustbox}

\numberwithin{equation}{section} 
\newtheorem{theorem}{Theorem}[section]
\newtheorem{corollary}[theorem]{Corollary}
\newtheorem{lemma}[theorem]{Lemma}
\newtheorem{proposition}[theorem]{Proposition}
\newtheorem{conjecture}[theorem]{Conjecture}
\newtheorem*{conjecture*}{Conjecture}
\theoremstyle{definition}
\newtheorem{definition}[theorem]{Definition}
\newtheorem*{definition*}{Definition}
\newtheorem{remark}[theorem]{Remark}
\newtheorem{example}[theorem]{Example}
\newtheorem{xxthm}{Theorem}

\newcommand{\cD}{\mathcal{D}}

\newcommand{\cF}{\mathcal{F}}

\newcommand{\cH}{\mathcal{H}}

\newcommand{\cK}{\mathcal{K}}
\newcommand{\cL}{\mathcal{L}}

\newcommand{\cO}{\mathcal{O}}
\newcommand{\cP}{\mathcal{P}}

\newcommand{\cR}{\mathcal{R}}
\newcommand{\cS}{\mathcal{S}}
\newcommand{\cT}{\mathcal{T}}
\newcommand{\cU}{\mathcal{U}}
\newcommand{\cV}{\mathcal{V}}
\newcommand{\cW}{\mathcal{W}}

\newcommand{\fI}{\mathfrak{I}}

\newcommand{\fU}{\mathfrak{U}}
\newcommand{\fV}{\mathfrak{V}}
\newcommand{\fW}{\mathfrak{W}}

\newcommand \sub{\mathrm{Sub}}
\newcommand \trans{\mathrm{Tr}}

\newcommand{\ind}{\mathrm{ind}}
\newcommand{\coind}{\mathrm{coind}}
\newcommand{\res}{\mathrm{res}}

\newcommand{\R}{\mathbb{R}}
\newcommand{\C}{\mathbb{C}}
\newcommand{\Z}{\mathbb{Z}}
\newcommand{\Irr}{\mathrm{Irr}}

\newcommand{\Hull}{\mathrm{Hull}}

\usepackage{tikz,tikz-cd}

\newcommand{\cpt}{
	\begin{tikzpicture}[scale=0.3,baseline=0.5mm]
		\node(0) at (0,0) {$\cdot$};
		\node(1) at (0,1) {$\cdot$};
	\end{tikzpicture}
}

\newcommand{\cpc}{
	\begin{tikzpicture}[scale=0.3,baseline=0.5mm]
		\node(0) at (0,0) {$\cdot$};
		\node(1) at (0,1) {$\cdot$};

        \draw (0,0) -- (0,1);  
	\end{tikzpicture}
}

\newcommand{\cppa}{
	\begin{tikzpicture}[scale=0.3,baseline=0.5mm]
		\node(0) at (0,0) {$\cdot$};
		\node(1) at (0,1) {$\cdot$};
		\node(2) at (0,2) {$\cdot$};
		
	\end{tikzpicture}
}

\newcommand{\cppb}{
	\begin{tikzpicture}[scale=0.3,baseline=0.2mm]
		\node(0) at (0,0) {$\cdot$};
		\node(1) at (0,1) {$\cdot$};
		\node(2) at (0,2) {$\cdot$};
		
		\draw (0,0) -- (0,1);
	\end{tikzpicture}
}

\newcommand{\cppc}{
	\begin{tikzpicture}[scale=0.3,baseline=0.5mm]
		\node(0) at (0,0) {$\cdot$};
		\node(1) at (0,1) {$\cdot$};
		\node(2) at (0,2) {$\cdot$};
		
		\draw (0,0) -- (0,1);
		\draw plot [smooth, tension=1.5] coordinates {(0,0) (0.5,1) (0,2)};
	\end{tikzpicture}
}

\newcommand{\cppd}{
	\begin{tikzpicture}[scale=0.3,baseline=0.5mm]
		\node(0) at (0,0) {$\cdot$};
		\node(1) at (0,1) {$\cdot$};
		\node(2) at (0,2) {$\cdot$};
		
		\draw (0,1) -- (0,2);
	\end{tikzpicture}
}

\newcommand{\cppe}{
	\begin{tikzpicture}[scale=0.3,baseline=0.5mm]
		\node(0) at (0,0) {$\cdot$};
		\node(1) at (0,1) {$\cdot$};
		\node(2) at (0,2) {$\cdot$};
		
		\draw (0,0) -- (0,1);
		\draw plot [smooth, tension=1.5] coordinates {(0,0) (0.5,1) (0,2)};
		\draw (0,1) -- (0,2);
	\end{tikzpicture}
}

\newcommand{\indt}{
	\begin{tikzpicture}[scale=0.2,baseline=0.2mm]
		\node(S) at (0,0) {$\cdot$};
		\node(W) at (-1,1) {$\cdot$};
		\node(E) at (1,1) {$\cdot$};
		\node(N) at (0,2) {$\cdot$};
		\path[-]
		
		;
	\end{tikzpicture}
}

\newcommand{\indpq}{
	\begin{tikzpicture}[scale=0.2,baseline=0.2mm]
		\node(S) at (0,0) {$\cdot$};
		\node(W) at (-1,1) {$\cdot$};
		\node(E) at (1,1) {$\cdot$};
		\node(N) at (0,2) {$\cdot$};
		\draw (0,0) -- (-1,1) ;
		\draw (0,0) -- (1,1) ;
	\end{tikzpicture}
}
\newcommand{\indpqp}{
	\begin{tikzpicture}[scale=0.2,baseline=0.2mm]
		\node(S) at (0,0) {$\cdot$};
		\node(W) at (-1,1) {$\cdot$};
		\node(E) at (1,1) {$\cdot$};
		\node(N) at (0,2) {$\cdot$};
		\draw (0,0) -- (1,1) ;
		\draw (-1,1) -- (0,2) ;
	\end{tikzpicture}
}
\newcommand{\indpqq}{
	\begin{tikzpicture}[scale=0.2,baseline=0.2mm]
		\node(S) at (0,0) {$\cdot$};
		\node(W) at (-1,1) {$\cdot$};
		\node(E) at (1,1) {$\cdot$};
		\node(N) at (0,2) {$\cdot$};
		\draw (0,0) -- (-1,1) ;
		\draw (1,1) -- (0,2) ;
	\end{tikzpicture}
}
\newcommand{\inds}{
	\begin{tikzpicture}[scale=0.2,baseline=0.2mm]
		\node(S) at (0,0) {$\cdot$};
		\node(W) at (-1,1) {$\cdot$};
		\node(E) at (1,1) {$\cdot$};
		\node(N) at (0,2) {$\cdot$};
		\draw (0,0) -- (-1,1) ;
		\draw (0,0) -- (1,1) ;
		\draw (0,0) -- (0,2) ;
		\draw (-1,1) -- (0,2) ;
		\draw (1,1) -- (0,2) ;
	\end{tikzpicture}
}
\newcommand{\indpqf}{
	\begin{tikzpicture}[scale=0.2,baseline=0.2mm]
		\node(S) at (0,0) {$\cdot$};
		\node(W) at (-1,1) {$\cdot$};
		\node(E) at (1,1) {$\cdot$};
		\node(N) at (0,2) {$\cdot$};
		\draw (0,0) -- (-1,1) ;
		\draw (0,0) -- (1,1) ;
		\draw (0,0) -- (0,2) ;
	\end{tikzpicture}
}
\newcommand{\indpqfp}{
	\begin{tikzpicture}[scale=0.2,baseline=0.2mm]
		\node(S) at (0,0) {$\cdot$};
		\node(W) at (-1,1) {$\cdot$};
		\node(E) at (1,1) {$\cdot$};
		\node(N) at (0,2) {$\cdot$};
		\draw (0,0) -- (-1,1) ;
		\draw (0,0) -- (1,1) ;
		\draw (0,0) -- (0,2) ;
		\draw (-1,1) -- (0,2) ;
	\end{tikzpicture}
}
\newcommand{\indpqfq}{
	\begin{tikzpicture}[scale=0.2,baseline=0.2mm]
		\node(S) at (0,0) {$\cdot$};
		\node(W) at (-1,1) {$\cdot$};
		\node(E) at (1,1) {$\cdot$};
		\node(N) at (0,2) {$\cdot$};
		\draw (0,0) -- (-1,1) ;
		\draw (0,0) -- (1,1) ;
		\draw (0,0) -- (0,2) ;
		\draw (1,1) -- (0,2) ;
	\end{tikzpicture}
}

\title[Equivariant linear isometries and infinite little discs operads]{Equivariant linear isometries and infinite little discs operads via transfer systems}

\begin{document}

\begin{abstract}
    In this article, we apply the recently developed theory of transfer systems to study the relationship between $G$-equivariant linear isometries and infinite little discs operads, for a finite group $G$. This framework allows us to reduce involved topological problems to discrete problems regarding the subgroup structure and representation theory of the group $G$. Our main result is an example of this: we classify the $G$-universes $\mathcal{U}$ for which the linear isometries operad $\mathcal{L}(\mathcal{U})$ and the infinite little discs operad $\mathcal{D}(\mathcal{U})$ are homotopically equivalent. To achieve this, we use ideas that originate from the work of Balchin-Barnes-Roitzheim on the combinatorics of transfer systems on a total order. Additionally, the use of transfer systems gives us insight into the algebraic structures that arise from equivariant homotopy theory. Compatible pairs of transfer systems provide rules for when multiplicative transfer maps can be paired with additive transfer maps. In the case that the group $G$ is abelian, we provide conditions for when the pair $(\mathcal{L}(\mathcal{U}),\mathcal{D}(\mathcal{U}))$ defines a maximally compatible pair of transfer systems. As a consequence, we contribute to a recent conjecture about equivariant operad pairs.    
\end{abstract}

\maketitle

\setcounter{tocdepth}{1}
\tableofcontents

\section*{Introduction}

\subsection*{Background} Homotopy coherent algebraic structures are abundant in modern algebraic topology. First defined by May in \cite{may1972geometry} to study the homotopy coherent algebra of iterated loop spaces, \textbf{operads} have served as a tool to encode these homotopy coherent structures. Let $\Sigma_n$ denote the $n$th symmetric group. An operad $\cO$ consists of, for each natural number $n$, a $\Sigma_n$-space of $n$-ary operations $\cO(n)$, along with additional data encoding identity operations, composition of operations and the subsequent unitality, associativity and $\Sigma_n$-equivariance conditions. The class of $E_\infty$-operads encode a topological generalisation of commutative monoids, where commutativity holds up to coherent homotopy. There are a few preferred $E_\infty$-operads, including the \textbf{infinite little discs operad} $\cD$ and the \textbf{linear isometries operad} $\cL$, which have classically appeared in the infinite loop space theory of \cite{may1977infinite}. While these particular operads may be used for their point-set properties, all $E_\infty$-operads are homotopically indistinguishable, see \cite{may1972geometry}. 

Group actions have had an increasing role in the development of homotopy theory. The most prominent example is the use of equivariant stable homotopy theory in the work of Hill-Hopkins-Ravenel on the Kervaire invariant problem, see \cite{hill2016nonexistence}; a surprising application to a question which, on its surface, had no connection to equivariant homotopy theory. $G$-operads, where one allows both $G$ and $\Sigma_n$ to act on the spaces $\cO(n)$, are useful tools to encode homotopy-coherent algebraic structures in $G$-spaces, see \cite{costenoble1991fixed}. Of course, one can import the $E_\infty$-operads $\cD$ and $\cL$ into the equivariant setting by requiring that $G$ acts trivially on the spaces $\cD(n)$ and $\cL(n)$. However, extra algebraic structure can arise by allowing $G$ to act non-trivially on the space of operations. Thus, it is preferable to also have less trivial equivariant generalisations of $\cD$ and $\cL$. To do this, one uses a certain type of infinite-dimensional $G$-inner-product space $\cU$, known as a \textbf{$G$-universe}, to define the equivariant infinite little discs operad $\cD(\cU)$ and the equivariant linear isometries operad $\cL(\cU)$, see \cite{guillou2017equivariant}. We wish to emphasise that there is a choice of universe $\cU$ and that the $G$-operads $\cD(\cU)$ and $\cL(\cU)$ may vary depending on this choice. Classically, one works over a \textbf{complete} $G$-universe $\cU_G$, where the operads $\cD(\cU_G)$ and $\cL(\cU_G)$ are homotopically equivalent. However, there can be universes for which $\cD(\cU)$ and $\cL(\cU)$ are distinct, see \cite[Theorem 4.22]{blumberg2015operadicmultiplicationsequivariantspectra}.

Motivated in part by the difference between $\cD(\cU)$ and $\cL(\cU)$, Blumberg and Hill defined $N_\infty$-operads in \cite{blumberg2015operadicmultiplicationsequivariantspectra} as an equivariant generalisation of $E_\infty$-operads. By \cite[Corollary 3.14]{blumberg2015operadicmultiplicationsequivariantspectra}, this class includes the equivariant infinite little discs $\cD(\cU)$ and linear isometries operads $\cL(\cU)$ over a universe $\cU$, which witness the fact that, unlike $E_\infty$-operads, there can be $N_\infty$-operads that are homotopically distinct. It has been shown in \cite{blumberg2015operadicmultiplicationsequivariantspectra, gutierrez2018encoding, bonventre2021genuine, Rubin_2021} that the homotopy category of $N_\infty$-operads is equivalent to a finite lattice $\mathrm{Ind}(G)$ of algebraic objects known as \textbf{indexing systems}. The lattice $\mathrm{Ind}(G)$ is further equivalent to the lattice $\trans(G)$ of \textbf{transfer systems}, as constructed by Rubin in \cite{rubin2021detecting}, and Balchin-Barnes-Roitzheim in \cite{balchin2022n}. Questions regarding the homotopy theory of $N_\infty$-operads can thus be reduced to discrete questions about either indexing systems or transfer systems. In  \cite[\S 4]{blumberg2015operadicmultiplicationsequivariantspectra}, Blumberg and Hill used this property of $N_\infty$-operads to establish that for all groups of order greater than three, there exist universes $\cU$ such that $\cL(\cU)$ and $\cD(\cU)$ are not equivalent, and that for all non-simple groups, there exists a universe $\cU$ such that $\cD(\cU)$ is not equivalent to the linear isometries operad $\cL(\cW)$ for any choice of $G$-universe $\cW$. They also gave important examples of where $\cD(\cU)$ and $\cL(\cU)$ are equivalent; if $\cU$ is isomorphic to the infinite regular representation $\R[G/N]^\infty$ of some quotient $G/N$ of $G$ by a normal subgroup $N$. In \cite{rubin2021detecting}, Rubin used the combinatorial nature of transfer systems to show that not every $N_\infty$-operad is equivalent to a linear isometries operad or an infinite little discs operad. He furthermore gave combinatorial properties of the transfer systems which arise from linear isometries and little discs operads. 

\subsection*{Results of this paper} Since the publications of \cite{blumberg2015operadicmultiplicationsequivariantspectra, rubin2021detecting}, the combinatorial theory of transfer systems has progressed significantly, see \cite{rubin2021operadic,franchere2022self,balchin2023lifting,bao2025transfer,adamyk2025minimalbaseshomotopicalcombinatorics}, and many of these developments have not yet been applied to the homotopy theory of $G$-operads. The present article has two aims. First, we wish to apply current knowledge about transfer systems to further understand the relationship between equivariant linear isometries and infinite little discs operads. Secondly, we aim to make contributions to the theory of transfer systems, with a view towards homotopical applications. We wish to emphasise that, despite our homotopical motivation, the main proofs of this article assume minimal prior knowledge of homotopy theory or operads. Instead, they rely on the subgroup structures and representation theory of finite groups.
With our aims in mind, the most obvious question is: For which universes $\cU$ are $\cL(\cU)$ and $\cD(\cU)$ equivalent? The main result of Section \ref{section:equivalence} is a complete classification of such universes.

\begin{xxthm}[Theorem \ref{thm:discEquivLinear}]
    \label{introthm:discEquivLinear}
    Let $G$ be a finite group and let $\cU$ be a $G$-universe. There exists a zig-zag of weak equivalences connecting $\cD(\cU)$ and $\cL(\cU)$ if and only if $\cU\cong \R[G/N]^\infty$ for some normal subgroup $N \unlhd G$.
\end{xxthm}

Note that $\cU \cong \R[G/N]^\infty$ is Blumberg and Hill's example of a universe where $\cD(\cU)$ and $\cL(\cU)$ are equivalent; the contribution of this article is a proof that this is the only example. In \ref{thm:bisaturated}, we also prove that if $\cD(\cU)$ and $\cL(\cW)$ are equivalent, for distinct universes $\cU$ and $\cW$, then they are also equivalent to $\cD(\R[G/N]^\infty)$ and $\cL(\R[G/N]^\infty)$ for some normal subgroup $N$. To prove these results, we completely classify the sublattice of bisaturated transfer systems, see Definition \ref{def:bisaturated}. This sublattice has independently appeared in a recent conjecture regarding the combinatorics of transfer systems in \cite{klanderman2025characterizingtransfersystemsnonabelian}.

Theorem \ref{introthm:discEquivLinear} raises another interesting question. If the linear isometries operad $\cL(\cU)$ is rarely equivalent to the little discs operad $\cD(\cU)$, then how else can we characterise the relationship between $\cL(\cU)$ and $\cD(\cU)$? One approach to such a question comes from the viewpoint of equivariant algebra. In \cite{blumberg2021}, Blumberg and Hill defined bi-incomplete Tambara functors; equivariant algebraic structures which generalise ordinary Tambara functors to settings where not all multiplicative and additive transfers are present. To control which multiplicative transfers can be paired with which additive transfers, Blumberg and Hill defined compatible pairs of indexing categories, which were translated to \textbf{compatible pairs} of transfer systems $(\cS,\cT)$ by Chan \cite{chan2024bi}. Here $\cS$ controls the multiplicative transfers and $\cT$ controls the additive transfers. In \cite{blumberg2021}, it was observed by Blumberg and Hill (in the equivalent context of indexing categories) that for a transfer system $\cT$, there is a unique maximal transfer system $\Hull(\cT)$ such that $(\Hull(\cT),\cT)$ form a compatible pair. The transfer system $\Hull(\cT)$ is known as the multiplicative hull of $\cT$, see \cite{blumberg2021,chan2025realizingcompatiblepairstransfer}. Let $\cT_\cO$ denote the transfer system corresponding to an $N_\infty$-operad $\cO$. In Section \ref{subsec:closureOperator}, we define a `closure operator' $\overline{(-)}$ on universes such that $\cD(\cU) \simeq \cD(\overline \cU)$, and, when the group $G$ is abelian, the multiplicative hull of $\cT_{\cD(\cU)}$ is equal to $\cT_{\cL(\overline\cU)}$.  
\begin{xxthm}[Theorem \ref{thm:maximalLinear}]
\label{introthm:maximalLinear}
    Let $G$ be an abelian group and let $\cU$ be a $G$-universe, then $\mathrm{Hull}(\cT_{\cD(\cU)})=\cT_{\cL(\overline\cU)}$.
\end{xxthm}
Theorem \ref{introthm:maximalLinear} is surprising at first, given how poor the combinatorial properties of linear isometries operads are, see \cite[Proposition 2.11]{rubin2021detecting}. After all, compatible pairs of transfer systems were not defined with the topology of equivariant operads in mind. However, Theorem \ref{introthm:maximalLinear} aligns with the recent program of connecting the purely algebraic notion of a compatible pair of transfer systems, as defined by Chan in \cite{chan2024bi}, with the more classical notion of operad pairs, originating from the work of May in \cite{may1977infinite}. The authors of \cite{chan2025realizingcompatiblepairstransfer} conjectured that every compatible pair of transfer systems can be suitably realised as an operad pair.  As an application of Theorem \ref{introthm:maximalLinear}, we have made the following contribution towards this conjecture.

\begin{xxthm}[Corollary \ref{cor:realizingCompatiblePairs}]
\label{introthm:realizingCompatiblePairs}
    Let $G$ be an abelian group, $\cT$ be a $G$ transfer system and $\cU$ be a $G$-universe. If $(\cT,\cT_{\cD(\cU)})$ forms a compatible pair of transfer systems, then there exists an $N_\infty$-operad pair $(\cP, \cO)$ with $\cT_\cP=\cT$ and $\cT_\cO = \cT_{\cD(\cU)}$.
\end{xxthm}

We conjecture that the assumption that $G$ is abelian can be removed from Theorems \ref{introthm:maximalLinear} and \ref{introthm:realizingCompatiblePairs}. Our proof involves working explicitly with the representation theory of abelian groups, using this to reduce to the cyclic quotients of $G$. The representation theory of non-abelian groups is significantly less well understood, thus we believe that a proof of the non-abelian case would require an entirely new strategy.  

Our remaining results address the combinatorics of linear isometries operads. Giving a combinatorial characterisation of the equivariant linear isometries has proven to be a rather difficult problem, see \cite{rubin2021detecting,hafeez2022saturatedlinearisometrictransfer,macbrough2023equivariant,xiaequivariant}. In \cite{blumberg2015operadicmultiplicationsequivariantspectra}, Blumberg and Hill noted that the indexing systems corresponding to equivariant linear isometries operads satisfy a sort of `horn-filling' property, which was translated into the notion of a saturated transfer system by Rubin in \cite{rubin2021detecting}. Rubin showed that among various groups of small order, there are some in which every saturated transfer system comes from a linear isometries operad, but that this need not hold in general. A group $G$ is saturated if every saturated $G$ transfer system comes from a linear isometries operad. Rubin conjectured that cyclic groups are saturated as long as the prime factors of their order are sufficiently large. This conjecture was proven in increasing generality in the sequence of papers \cite{rubin2021detecting,hafeez2022saturatedlinearisometrictransfer,macbrough2023equivariant}. However, the question of which non-abelian groups are saturated has not been as well studied, apart from the work of Xia in \cite{xiaequivariant}, who proved that certain non-abelian $p$-groups are not saturated. We generalise their strategy and prove that certain infinite families of non-abelian groups are not saturated. 

\begin{xxthm}[Corollary \ref{cor:symmetricNotsaturated}]
    The symmetric group $\Sigma_n$ is saturated if and only if $n \leq 2$, and the alternating group $A_m$ is saturated if and only if $n \leq 3$.
\end{xxthm}

\begin{xxthm}[Corollary \ref{cor:hamiltonianNotSaturated}]
    Let $G$ be a finite Hamiltonian group, that is, $G$ is non-abelian, and every subgroup of $G$ is normal, then $G$ is not saturated.
\end{xxthm}

These theorems are applications of the results of Section \ref{section:saturated}, which could be used to deduce that additional families of non-abelian groups are not saturated. In particular, no examples of saturated non-abelian groups have been found. We conjecture that this is because no such groups exist. 

\begin{conjecture*}
    If $G$ is a saturated group, then $G$ is abelian.
\end{conjecture*}

If the above conjecture is true, then a saturated group must be abelian of rank at most 2, where only the rank 2 case remains open in general, see \cite{macbrough2023equivariant} for a proof of this, and a discussion of the rank 2 case.

\subsection*{Notation and Conventions} Throughout the present paper, all groups will be finite, usually denoted $G$, with subgroups usually denoted $L, K, H$. For two groups $H,G$, we denote by $H \leq G$ that $H$ is a subgroup of $G$. We denote by $H\unlhd G$ that $H$ is a normal subgroup of $G$. The lattice of subgroups is denoted by $\sub(G)$, and the lattice of normal subgroups is denoted by $\sub_\unlhd(G)$. If $H \leq G$ and $g \in G$, then we use the notation $H^g := gHg^{-1}$ for conjugate subgroups. If also $K \leq G$, then we use the notation $KH := \{kh \mid k\in K, h \in H\}$ and  $K \cap H :=\{g \in G \mid g \in K \text{ and } g\in H\}$. 

We will use fonts to differentiate between finite-dimensional $G$-representations $U, V, W$, (real) $G$-universes $\cU, \cV, \cW$ and complex $G$-universes $\fU, \fV, \fW$. All representations and universes are assumed to be real unless explicitly stated otherwise. If $\mathbb{K}$ is a field, either $\R$ or $\C$ in this paper, then $\Irr_\mathbb{K}(G)$ denotes the set of isomorphism classes of irreducible $G$-representations over $\mathbb{K}$. If $n \in \mathbb{N}$, then $\Irr_\mathbb{K}^n(G) \subseteq \Irr_\mathbb{K}(G)$ is the subset of $n$-dimensional irreducible representations and $\Irr_\mathbb{K}^{>n}(G) \subseteq \Irr_\mathbb{K}(G)$ is the subset of representations of dimension greater than $n$. If $\mathbb{K}$ is dropped in the above notation, then it is assumed that $\mathbb{K}=\R$.

Unless stated otherwise, we will use $\cO,\cP$ to denote $N_\infty$-operads. We use $\cT,\cS,\cR$ to denote transfer systems with order relation $\to_\cT$, $\to_\cS$,$\to_\cR$. We use $\cT_\cO$ to denote the transfer system associated to an $N_\infty$-operad $\cO$, whose order relation is simply denoted $\to_\cO$.

\subsection*{Acknowledgements} The present article was written as part of the author's PhD research, under the supervision of Magdalena K\c{e}dziorek. The author would like to thank Magdalena K\c{e}dziorek and Niall Taggart for their guidance and for many helpful conversations. The author also thanks the referee for a swift report and several useful suggestions. This work was partially supported by the Dutch Research Council (NWO) grant number OCENW.M.22.358.

\section{\texorpdfstring{$N_\infty$}--operads and transfer systems}

We begin with a brief recollection of the theory of $N_\infty$-operads and transfer systems, following \cite{blumberg2015operadicmultiplicationsequivariantspectra, rubin2021detecting, rubin2021operadic, balchin2024combinatorics}. Along the way, we will prove some prerequisite results which do not currently appear in the literature.

\subsection{$N_\infty$-operads} 

\begin{definition}[{\cite[Definition 2.1]{costenoble1991fixed}}]
    \label{def:gOperad}
    A \textbf{$G$-operad} $\cO$ consists of a $G \times \Sigma_n$-space $\cO(n)$ for each integer $n\geq 0$, together with: 
    \begin{enumerate}
        \item A $G$-fixed identity element $1 \in \cO(1)$.
        \item For each $k \geq 0$ and $n_1,\dots ,n_k \geq 0$, a $G$-equivariant composition map 
        \[
        \cO(k) \times \cO(n_1) \times \cdots \times \cO(n_k) \to \cO(n_1 + \cdots + n_k).
        \]
    \end{enumerate}
    These data are subject to compatibility conditions encoding unitality, associativity, and permutation equivariance.
    If $\cO$ and $\cP$ are two $G$-operads, a $G$-operad map $\cO \to \cP$ consists of $G \times \Sigma_n$-equivariant continuous maps $\cO(n) \to \cP(n)$, which are compatible with the composition maps and identity element. Let $\mathrm{Op}^G$ denote the category whose objects are $G$-operads and morphisms are $G$-operad maps.
\end{definition}

For a homotopical analysis of the category $\mathrm{Op}^G$, we lift the notion of weak equivalence from spaces. 
\begin{definition}[{\cite[Definition 3.9]{blumberg2015operadicmultiplicationsequivariantspectra}}]
    \label{def:weakEquivalence}
    A map $\cO \to \cP$ of $G$-operads is a \textbf{weak equivalence} if each map $\cO(n)^\Gamma \to \cP(n)^\Gamma$ is an equivalence for all subgroups $\Gamma \leq G \times \Sigma_n$.
\end{definition}

The present paper is concerned with two particular classes of $G$-operads which have played an important role in the development of equivariant stable homotopy theory. Namely, the equivariant infinite little discs and linear isometries operads, which provide point-set models for additive and multiplicative structures in $G$-spectra.   

\begin{definition}
    \label{def:littleDiscs}
    Let $V$ be a finite-dimensional orthogonal $G$-representation. Let $D(V)$ denote the unit disc of $V$ centred at the origin. A \textbf{little $V$-disc} in $D(V)$ is an affine, but not necessarily equivariant, map $D(V) \to D(V)$. The \textbf{little $V$-discs operad} $\cD(V)$ is defined in arity $n$ as the space of $n$ disjoint little $V$-discs in $D(V)$, with $G$-action by conjugation. 
\end{definition}

Note that if $V \hookrightarrow W$ is an inclusion of representations, there is an induced map of operads $\cD(V) \to \cD(W)$ given by `thickening' little $V$-discs into little $W$-discs, see \cite[Example 2.2(c)]{costenoble1991fixed}.
In the non-equivariant context, the infinite little discs operad $\cD$ is defined as a colimit of the little $n$-discs operads (take $V = \R^n$ in Definition \ref{def:littleDiscs}) taken over all $n \geq 0$. The equivariant analogue of this definition requires taking a colimit over a suitable infinite collection of finite-dimensional representations. Such data are typically packaged into a universe.

\begin{definition}
    \label{def:universe}
    A \textbf{$G$-universe} $\cU$ is a countably infinite-dimensional orthogonal $G$-representation, such that $\cU$ contains a non-zero subrepresentation on which $G$ acts trivially, and for every finite-dimensional subrepresentation $V \subset \cU$, the $G$-representations $V^{\oplus n}$, with $n \in \mathbb{N}$, are also subrepresentations of $\cU$. A $G$-universe $\cU$ is \textbf{complete} if every finite dimensional $G$-representation embeds into $\cU$. A $G$-universe $\cU$ is \textbf{trivial} if it has the trivial $G$-action. 
\end{definition}

In particular, the infinite regular representation $\R[G]^\infty$ is a complete universe and $\R^\infty$ is a trivial $G$-universe. More generally, every $G$-universe is isomorphic to one of the form $(\R \oplus \bigoplus_{i=0}^n V_i)^\infty$, where each $V_i$ is an irreducible $G$-representation.  

\begin{definition}
    For a universe $\cU$, the \textbf{infinite little discs operad} $\cD(\cU)$ is defined as the colimit of $\cD(V)$ taken over the poset of finite-dimensional subrepresentations $V \subseteq \cU$. 
\end{definition}

\begin{definition}
    For a universe $\cU$, the \textbf{linear isometries operad} $\cL(\cU)$ is defined in arity $n$ as the space of (not necessarily equivariant) linear isometries $\cU^{\oplus n} \to \cU$, with $G$-action by conjugation.
\end{definition}

Rather than work in the category of $G$-operads, we shall restrict ourselves to the much more convenient subcategory of $N_\infty$-operads, which contains all of the operads we are interested in for the present paper. 

\begin{definition}
    \label{def:family}
    A \textbf{family} of subgroups is a collection $\cF$ of subgroups of $G$, which is closed under conjugation and passage to subgroups. 
    A universal space for a family $\cF$ is a $G$-space $E \cF$ such that: 
    $$ (E \cF)^H \simeq \begin{cases}
        \ast,  & \text{if } H \in \cF, \\
        \emptyset, & \text{otherwise.}
    \end{cases} $$
\end{definition}

\begin{definition}[{\cite[Definition 3.7]{blumberg2015operadicmultiplicationsequivariantspectra}}]
\label{def:Ninfty}
    An \textbf{$N_\infty$-operad} is a $G$-operad $\cO$ such that: 
    \begin{enumerate}
        \item For all $n \geq 0$, the $\Sigma_n$ action on $\cO(n)$ is free. 
        \item For all $n \geq 0$, the space $\cO(n)$ is a universal space for a family of subgroups of $G \times \Sigma_n$ containing each subgroup of the form $H \times \{e\}$ where $H$ is a subgroup of $G$.
        \item $\cO(0)$ is $G$-contractible. 
    \end{enumerate}
    Let $N_\infty^G$ denote the full subcategory of $\mathrm{Op}^G$ consisting of $G$-$N_\infty$-operads. Denote by $\mathbf{Ho}(N_\infty^G)$ the homotopy category obtained by localising $N_\infty^G$ with respect to the weak equivalences in Definition $\ref{def:weakEquivalence}$.
\end{definition}

Under certain conditions, see \cite[\S7]{blumberg2015operadicmultiplicationsequivariantspectra}, an $N_\infty$-operad $\cO$ induces `wrong way' maps $X^K \to X^H$ between the fixed points of an $\cO$-algebra $X$ for $K \leq H \leq G$. The pairs $(K,H)$, for which such wrong-way maps are induced, are invariant under the weak equivalences of Definition \ref{def:weakEquivalence}. Thus, one can use this data as a combinatorial invariant for $N_\infty$-operads. This collection of pairs of subgroups can be packaged into a partial order $\to_\cO$ on the set of subgroups $\sub(G)$ of $G$, whose properties are captured by the definition of a transfer system.

\begin{definition}[{\cite{rubin2021detecting, balchin2022n}}]
    \label{def:transferSystem}
    Let $G$ be a finite group. A \textbf{$G$ transfer system} $\cT$ consists of a partial order $\to_\cT$ on the set of subgroups $\sub(G)$ of $G$ such that:
    
    \begin{enumerate}
        \item $\cT$ refines inclusions; if $K \to_\cT H$ then $K \leq H$.  
        \item $\cT$ is closed under restriction; if $L \to_\cT H$ and $K \leq H$ then $(K\cap L) \to_\cT K$
        \item $\cT$ is closed under conjugation; if $K \to_\cT H$ and $g \in G$, then $K^g \to_\cT H^g$.
    \end{enumerate}
    We will usually call a relation $K \to_\cT H$ a `transfer' in $\cT$.
    We denote by $\trans(G)$ the set of $G$ transfer systems.
    \end{definition}

    \begin{example}[Transfer system diagrams for cyclic groups]
    \label{ex:diagrams}
    Since transfer systems are posets, they can be represented by Hasse-diagrams. In particular, these diagrams are subgraphs of the Hasse-diagram for the lattice of subgroups. We include examples below for the cyclic groups $C_n$ of order $n$. 

    \[\begin{tabular}{c|c|c}
    $G$ & $\mathrm{Sub}(G)$ & $\trans(G)$ \\
    \hline
    \hline
       $C_2$ & 
       \begin{tikzcd}[row sep=small, column sep=small]
	    {C_2} \\
	    e
	    \arrow[from=2-1, to=1-1]
        \end{tikzcd} & $\cpt \quad \cpc$ \\
    \hline 
    $C_4$ & \begin{tikzcd}[row sep=small, column sep=small]
	{C_4} \\
	{C_2} \\
	e
	\arrow[from=2-1, to=1-1]
	\arrow[curve={height=12pt}, from=3-1, to=1-1]
	\arrow[from=3-1, to=2-1]
\end{tikzcd} & $\cppa \quad \cppb \quad \cppc \quad \cppd \quad \cppe$ \\
    \hline
    $C_6$ & \begin{tikzcd}[row sep=small, column sep=tiny]
	& {C_6} & \\
	{C_3} && {C_2} \\
	& e
	\arrow[from=2-1, to=1-2]
	\arrow[from=2-3, to=1-2]
	\arrow[from=3-2, to=1-2]
	\arrow[from=3-2, to=2-1]
	\arrow[from=3-2, to=2-3]
\end{tikzcd} & $\indt \quad \indpq \quad \indpqf \quad \indpqp \quad \indpqq \quad \indpqfp \quad \indpqfq  \quad \inds$
  
    \end{tabular}\]
    \end{example}

    One can define a partial order $\leq$ on the set $\trans(G)$ such that $\cT\leq \cS$ whenever $\cT$ refines $\cS$. By `refines' we mean that whenever $K \to_{\cT} H$, we also have that $K \to_{\cS} H$. We may now view the poset $(\trans(G),\leq)$ as a category, which we shall simply denote as $\trans(G)$ from now on. 
    
    \begin{theorem}[\cite{blumberg2015operadicmultiplicationsequivariantspectra,gutierrez2018encoding,bonventre2021genuine,Rubin_2021}]
    \label{thm:homotopyCategory}
        There is an equivalence of categories:
        \[
        \mathbf{Ho}(N_\infty^G) \to \trans(G),
        \]
        from the homotopy category of $N_\infty$-operads to the lattice of transfer systems. 
    \end{theorem}

    \begin{remark}
        The above theorem was originally stated in a different, but equivalent, form using the poset of indexing systems $\mathrm{Ind}(G)$. It was later shown that there exists an isomorphism $\mathrm{Ind}(G) \to \trans(G)$ \cite[Theorem 3.7]{rubin2021detecting}.
    \end{remark}

    We will use $\cT_\cO$ to denote the transfer system corresponding to (the weak homotopy type of) an $N_\infty$-operad $\cO$, and we shall write $\to_\cO = \to_{\cT_\cO}$ as shorthand. By Theorem \ref{thm:homotopyCategory}, $\cT_\cO \leq \cT_\cP$ may be interpreted operadically as ``there exists a unique map $\cO \to \cP$ in the homotopy category of $G$-operads" and $\cT_\cO = \cT_\cP$ may be interpreted as ``there exists a zig-zag of weak equivalences connecting $\cO$ and $\cP$".

Blumberg and Hill \cite[\S 4]{blumberg2015operadicmultiplicationsequivariantspectra} observed that the functor $\mathbf{Ho}(N^G_\infty)\to \mathrm{Ind(G)}$, which was only known to be fully faithful at the time, could be used to examine the relationship between different point-set models of $N_\infty$-operads, rather than constructing explicit maps between $G$-operads. In particular, they established characterisations of the indexing systems corresponding to linear isometries and infinite little discs operads, which translate to the following conditions when passing to transfer systems. 

\begin{proposition}[{\cite[Theorems 4.18 \& 4.19]{blumberg2015operadicmultiplicationsequivariantspectra}}]
    \label{prop:lineardiscTransfersystem}
    Let $\cU$ be a $G$-universe, and let $K \leq H \leq G$ be subgroups of $G$, then: 
    \begin{itemize}
        \item $K \to_{\cD(\cU)} H$ if and only if $H/K$ embeds $H$-equivariantly into $\res^G_H\cU$.
        \item $K \to_{\cL(\cU)} H$ if and only if $\ind^H_K \res^G_K \cU$ embeds $H$-equivariantly into $\res^G_H\cU$.
    \end{itemize}
\end{proposition}

All of the following examples were already done in \cite{rubin2021detecting}, where the reader can find several additional examples. It is a useful exercise to try to classify these transfer systems for various small groups. 

\begin{example}[$G=C_2$]
    There are two irreducible real $C_2$-representations; the one-dimensional trivial representation $\R$ and the one-dimensional sign representation $\sigma$. Since every universe must contain a trivial summand, there are two universes up to isomorphism; $\R^\infty$ and $(\R \oplus \sigma)^\infty$. The only orbits which embed into $\R$ are the trivial orbits, and all orbits embed into $\sigma$. Furthermore, $\ind^{C_2}_{\{e\}}\R$ is the regular representation $\R[C_2] \cong \R\oplus \sigma$. We can then produce the following table: 
    
    \[\begin{tabular}{c|c|c}
        $\cU$ & $\cT_{\cD(\cU)}$ & $\cT_{\cL(\cU)}$ \\
        \hline
        \hline
        $\R^\infty$ &  $\cpt$  & $\cpt$ \\
        \hline
        $(\R \oplus \sigma)^\infty$ & $\cpc$ & $\cpc$
    \end{tabular}\]
    
\end{example}

\begin{example}[$G=C_4$]
     \label{c4linearsteiner}
     There are three irreducible $C_4$-representations: the one-dimensional trivial representation $\R$, the one-dimensional sign representation $\sigma$ on $C_4/C_2 \cong C_2$ and the two-dimensional representation $\lambda_4$ where the generator acts by a rotation of order $4$. The only orbit which embeds into $\sigma$ is $C_4/C_2$, and $\lambda_4\setminus \{0\}$ has a free $C_4$ action, so the only non-trivial orbits which embed into $\lambda_4$ are $C_4/\{e\}$ and $C_2/\{e\}$. Moreover, $\ind^{C_4}_{C_2}\res^{C_4}_{C_2}\sigma \cong \ind^{C_4}_{C_2}\R \cong \R[C_4/C_2]$ and $\ind^{C_4}_{C_2} \res^{C_4}_{C_2}\lambda_4\cong\ind^{C_4}_{C_2}(\sigma \oplus \sigma)\cong \lambda_4\oplus \lambda_4$.
     
    \[\begin{tabular}{c|c|c}
        $\cU$ & $\cT_{\cD(\cU)}$ & $\cT_{\cL(\cU)}$ \\
        \hline
        \hline
        $\R^\infty$ &  $\cppa$  & $\cppa$ \\
        \hline
        $(\R \oplus \sigma)^\infty$ & $\cppd$ & $\cppd$ \\
        \hline
        $(\R \oplus \lambda_4)^\infty$ & $\cppc$ & $\cppb$ \\
        \hline
        $(\R \oplus \sigma \oplus \lambda_4)^\infty$ & $\cppe$ & $\cppe$
    \end{tabular}\]
\end{example}

\begin{example}[$G=C_6$]
    \label{c6linearsteiner}
    There are four irreducible $C_6$ representations: the trivial one-dimensional representation $\R$, the one-dimensional $C_6/C_3 \cong C_2$ sign representation $\sigma$, and the two-dimensional representations $\lambda_3$ and $\lambda_6$ where the generator acts by rotations of order $3$ and $6$ respectively. 
    \[\begin{tabular}{c|c|c}
        $\cU$ & $\cT_{\cD(\cU)}$ & $\cT_{\cL(\cU)}$\\
        \hline
        \hline
        $\R^\infty$ &  $\indt$  & $\indt$ \\
        \hline
        $(\R \oplus \sigma)^\infty$ & $\indpqq$ & $\indpqq$ \\
        \hline
        $(\R \oplus \lambda_3)^\infty$ & $\indpqp$ & $\indpqp$ \\
        \hline
        $(\R \oplus \lambda_6)^\infty$ & $\indpqf$ & $\indpq$ \\
        \hline
        $(\R \oplus \sigma \oplus \lambda_3)^\infty$ & $\inds$ & $\indpq$ \\
         \hline
        $(\R \oplus \sigma \oplus \lambda_6)^\infty$ & $\indpqfq$ & $\indpq$ \\
        \hline
        $(\R \oplus \lambda_3 \oplus \lambda_6)^\infty$ & $\indpqfp$ & $\indpq$
        \\
        \hline
        $(\R \oplus \sigma \oplus\lambda_3 \oplus \lambda_6)^\infty$ & $\inds$ & $\inds$
        
    \end{tabular}\]
     
\end{example}

The reader may have noticed in the above examples that the condition $K \to_{\cL(\cU)} H$ implies that $K \to_{\cD(\cU)} H$. One can verify that this is true for all groups $G$. Thus, we can deduce the following theorem.  

\begin{theorem}[{\cite[Corollary 4.20]{blumberg2015operadicmultiplicationsequivariantspectra}}]
    \label{thm:discsContainsLinear}
    For each $G$-universe $\cU$, there exists a unique map $\cL(\cU) \to \cD(\cU)$ in the homotopy category of $G$-operads.
\end{theorem}

The above theorem illustrates the utility of transfer systems as a tool to study the homotopy theory of $N_\infty$-operads, which we will exploit in order to prove results in the spirit of Theorem \ref{thm:discsContainsLinear}.

\subsection{Minimal fibrant subgroups} We will now discuss some useful notions and results about transfer systems.

\begin{proposition}
    If $G$ is a finite group and $\cT$ is a $G$ transfer system, then there exists a unique minimal subgroup $M$ such that $M \to_\cT G$. The subgroup $M$ is known as the minimal fibrant subgroup of $\cT$. 
\end{proposition}
\begin{proof}
    An analogous proof is given in \cite[Lemma 2.2]{balchin2024combinatorics} in the setting of transfer systems on a lattice. Suppose $H$ and $K$ are subgroups of $G$ such that $H \to_\cT G$ and $K \to_\cT G$. By restriction of $K \to_\cT G$ to $H$, we have $H \cap K \to_\cT H$. Thus, by transitivity, we obtain the transfer $H \cap K \to_\cT G$. Since $G$ has finitely many subgroups, one obtains $M$ as the intersection of all subgroups $H$ with the property that $H \to_\cT G$.
\end{proof}

The notion of a minimal fibrant subgroup is crucial to many of the existing counts of transfer systems, see \cite{balchin2022n, balchin2024combinatorics, bao2025transfer}. The minimality of $M$, along with the closure under restriction and transitivity of the transfer system $\cT$, ensures that there can be no transfers $K \to_\cT H$ with $M \nleq K $ and $M \leq H$. This effectively splits the transfer system $\cT$ into a kind of disjoint sum. In the case that $G= C_{p^n}$ for some prime number $p$ and some natural number $n$, this observation was used in \cite[Theorem 20]{balchin2022n} to establish the reccurence relation 
\[
|\trans(C_{p^n})|=\sum_{i+j = n-1}\left(|\trans(C_{p^i})| + |\trans(C_{p^j})|\right),
\]
and, in turn, show that $|\trans(C_{p^n})|$ is the $(n+1)$th Catalan number. The applications of minimal fibrant subgroups in \cite{balchin2022n, balchin2024combinatorics, bao2025transfer} have either been in the abelian setting, in \cite{balchin2022n, balchin2024combinatorics, bao2025transfer} or have lifted transfer systems over non-abelian groups from transfer systems on a lattice in \cite{balchin2024combinatorics}, using the ideas of \cite{balchin2023lifting}. In the present paper, we also wish to examine general non-abelian groups, where one obtains the following restriction.

\begin{proposition}
    The minimal fibrant subgroup of a transfer system $\cT$ is always a normal subgroup of $G$. In particular, taking the minimal fibrant subgroup of a transfer system yields an order-reversing map 
    \[
    \trans(G) \to \sub_\unlhd(G),
    \]
    where we use $\sub_\unlhd(G)$ to denote the lattice of normal subgroups of $G$.
\end{proposition}
\begin{proof}
    Suppose that $H$ is a subgroup such that $H \to_\cT G$, and pick any $g \in G$. By closure under restriction, we see that $H \cap H^g\to_\cT G$. If $H^g \neq H$, then $H \cap H^g$ is strictly contained in $H$. Consequently, if $M$ is minimal with the property that $M \to_\cT G$, then $M$ must be normal. 

    To see that the map $\trans(G) \to \sub_\unlhd(G)$ is order reversing, recall that an order relation $\cT \leq \cS$, in $\trans(G)$, means precisely that $K \to_{\cT} H$ implies $K \to_{\cS} H$. Thus, if $H \to_{\cT} G$, we have that $H \to_{\cS} G$. It follows that the minimal fibrant subgroup of $\cS$ is contained in the minimal fibrant subgroup of $\cT$.
    \end{proof}
    
\subsection{Saturated and disc-like transfer systems}
 Rubin \cite{rubin2021detecting} showed that the transfer systems corresponding to linear isometries operads and the little discs operads belong to particular classes of transfer systems known as saturated and disc-like (or cosaturated) transfer systems, respectively.   

\begin{definition}
    \label{def:saturatedTransfer}
    A transfer system $\cT$ is \textbf{saturated} if, whenever $L \leq K \leq H \leq G$, and $L \to_\cT H$, then $K \to_\cT H$.
\end{definition}

The additional condition that saturated transfer systems satisfy can be seen as a two-out-three property, in the sense that whenever a transfer system $\cT$ is saturated, then the existence of two of the three transfers $L \to_\cT K$, $K \to_\cT H$, $L \to_\cT H$ implies the existence of the third. 

\begin{definition}
    \label{def:disclikeTransfer}
    A transfer system $\cT$ is \textbf{disc-like} (or \textbf{cosaturated}) if it is generated by a set of transfers $\{H_i \to G\}_{i \in I}$. 
\end{definition}

\begin{remark}
    By `generated', we mean the unique smallest transfer system containing a given binary relation $\cR$, which we denote by $\langle \cR \rangle$. The notion of a generated transfer system was first used by Rubin in order to describe the behaviour of infinite little discs operads. An algorithm for computing the transfer system generated by a given relation is given in \cite[Construction A.1]{rubin2021detecting}, which has since been applied to the combinatorics of transfer systems in \cite{adamyk2025minimalbaseshomotopicalcombinatorics}.
\end{remark}

\begin{theorem}[\cite{rubin2021detecting}]
    For a $G$-universe $\cU$, the transfer system $\cL(\cU)$ is saturated and the transfer system $\cD(\cU)$ is disc-like.
\end{theorem}

Saturated and disc-like transfer systems have many desirable properties. As such, these classifications make the combinatorial study of linear isometries and little discs operads significantly more manageable. The propositions below are examples of such desirable properties, which will be required in Section \ref{subsec:changeGroup}.

\begin{proposition}
    \label{prop:disclikeRestriction}
    If $G$ is a finite group, and $\cT$ is a disc-like transfer system over $G$, then $\cT$ can be obtained by closing the following set of subgroup inclusions under restriction: 
    \[
    S=\{H \to G \mid H \in \mathrm{Sub}(G) \text{ and }H \to_\cT G\}.
    \]   
\end{proposition}
\begin{proof}
    Since $\cT$ is known to be disc-like, it is generated by a set of subgroup inclusions with target $G$; since any such set is contained in $S$ by construction, we can conclude that $\to_\cT$ is generated by $S$.
    Rubin's algorithm \cite[Theorem A.2]{Rubin_2021} states that $\cT$ can be obtained by closing $S$ under conjugation, restriction and transitivity, in that order. By construction, $S$ is already closed under conjugation. Thus, it suffices to show that its closure under restriction is closed under transitivity.

    Indeed, suppose that $A\to B$ and $B \to C$ are both in the closure of $S$ under restriction. That is, there exist subgroups $K, H \leq G$ such that $K \to_\cT G$, $H \to _\cT G$, $H \cap C = B$ and $K \cap B = A$. By replacing $K$ with $K \cap H$, we may assume that $K \leq H$. It follows that 
    $
    K \cap C= (K\cap H) \cap C = K \cap B = A
    $
    as desired. 
\end{proof}

\begin{proposition}
    \label{prop:saturatedNormalCore}
    For a finite group $G$, and a subgroup $H$, let \[
    \mathrm{Core}_G(H)=\bigcap_{g \in G} H^g
    \] denote the normal core of $H$ in $G$. If $\cT$ is a saturated $G$ transfer system, then $K \to_\cT H$ if and only if $\mathrm{Core}_H(K) \to_\cT H$.
\end{proposition}
\begin{proof}
    The `if' part of the proposition immediately follows from the saturated property of $\cT$. The `only if' part holds for all transfer systems $\cT$. Indeed, closure under conjugacy ensures that $K^h \to_\cT H$ for all $h \in H$. By inductively applying closure under restriction and transitivity, we see that $\mathrm{Core}_H(K) \to_\cT H$, as desired.
\end{proof}

\subsection{The lattice of universes}

Given a finite group $G$, one can form a finite lattice $\mathbf{Uni}_\R(G)$ consisting of isomorphism classes of real $G$-universes, which we will denote by $\mathbf{Uni}(G)$ unless it needs to be distinguished from the lattice of complex universes that will be defined shortly. We have that $[\cU] \leq [\cW]$ if there exists an embedding $\cU \hookrightarrow \cW$. Joins are given by direct sums, and the meet $[\cU] \wedge [\cW]$ is the universe generated by all irreducible $G$-representations which embed into both $\cU$ and $\cW$. It is natural to assume that the assignments $\cU \mapsto \cL(\cU)$ and $\cU \mapsto \cD(\cU)$ would induce a pair of lattice maps $\cT_\cL,\cT_\cD : \mathbf{Uni}(G) \to \trans(G)$. However, it was established by Rubin that this is not true in general.

\begin{proposition}[{\cite[Proposition 2.10]{rubin2021detecting}}]
    The map $\cT_\cD : \mathbf{Uni}(G) \to \trans(G)$ preserves the order, joins, the top element and the bottom element, but does not preserve meets in general. 
\end{proposition}

\begin{proposition}[{\cite[Proposition 2.11]{rubin2021detecting}}]
    The map $\cT_\cL : \mathbf{Uni}(G) \to \trans(G)$ preserves the top element and the bottom element but does not preserve the order, meets or joins in general.
\end{proposition}

It is often more convenient to work with representations over an algebraically closed field. For this reason, one can define the notion of a complex $G$-universe. That is, a countably infinite-dimensional complex $G$-representation into which every finite-dimensional subrepresentation embeds countably infinitely many times. Analogously to the real case, one can define a lattice $\mathbf{Uni}_\C(G)$ of isomorphism classes of complex $G$-universes. A treatment of the lattice of complex universes as the subset lattice of $\Irr_\C(G)$ can be found in \cite{macbrough2023equivariant}. We shall sketch the corresponding geometric picture below.

There is a natural forgetful meet-preserving map $F : \mathbf{Uni}_\C(G) \to \mathbf{Uni}_\R(G)$, by restricting scalars along the usual inclusion $\R \hookrightarrow \C$. Given a complex $G$-universe $\fU$, one can define the infinite little discs operad as follows $\cD(\fU) := \cD(F(\fU))$, such that $K \to_{\cD(\fU)} H$ if and only if there is an $H$-equivariant embedding $H/K \hookrightarrow \fU$. However, to define a suitable $\cL(\fU)$, one has to be more careful. This is because linear maps between complex vector spaces do not correspond to linear maps between the underlying real vector spaces. For the present article, it is sufficient to define the transfer system $\cT_{\cL(\fU)}$ by combinatorial means. Specifically, we define the order relation $\to_{\cL(\fU)}$, such that $K \to_{\cL(\fU)} H$ if and only if $\ind^H_K\res^G_K\fU$ $H$-equivariantly embeds into $\res^G_H \fU$. Here $\ind$ and $\res$ denote the induction and restriction functors between complex representations, and $\to_{\cL(\fU)}$ defines a transfer system ${\cL(\fU)}$, analogously to the real case. We have an injective complexification map $(-) \otimes \C:\mathbf{Uni}_\R(G) \to \mathbf{Uni}_\C(G)$ which commutes with induction and restriction. Hence, for a real $G$-universe $\cU$, we have $\cT_{\cL(\cU)} = \cT_{\cL(\cU \otimes \C)}$ and the isomorphism class of $\cU$ is uniquely determined by the isomorphism class of $\cU \otimes \C$.

\subsection{Change-of-group maps}
\label{subsec:changeGroup}
In this section, we will give an outline of the change-of-group functors which were constructed in \cite[\S 5]{rubin2021operadic}. For the present article, we will specialise to the cases of subgroup inclusions and quotient maps, where we give direct definitions corresponding to Rubin's constructions.  
Given a group homomorphism $f: G \to G'$, Rubin has defined maps $f_L,f_R : \trans(G) \to \trans(G')$ and $f_L^{-1},f_R^{-1} : \trans(G') \to \trans(G)$, such that we have adjunctions $f_L \dashv f_R^{-1}$ and $f_L^{-1} \dashv f_R$, see \cite[\S 5]{rubin2021operadic}. By \cite[Proposition 5.13]{rubin2021operadic}, there is an equality $f^{-1}_L = f^{-1}_R$ if and only if $f$ is injective. In particular, if $\iota: H \to G$ is an inclusion of a subgroup $H$. Then we call $\iota^{-1}_L = \iota^{-1}_R$ restriction to $H$ and denote this as $\res^G_H$. We call the left and right adjoints of $\res^G_H$ induction and coinduction, denoted by $\ind^G_H$ and $\coind^G_H$ respectively. See \cite[\S 6]{rubin2021operadic} for a discussion of how $\ind^G_H$, $\res^G_H$ and $\coind^G_H$ lift to functors on the level of operads. The following characterisations of $\res^G_H$ and $\ind^G_H$ can be deduced from \cite[Lemma 5.7]{rubin2021operadic}.

\begin{proposition}
    \label{prop:resindCharacterisation}
    If $\cT$ is a $G$ transfer system, then $\res^G_H\cT$ is defined by the rule that $L\to_{\res^G_H\cT}K$ if and only if $L\leq K \leq H$ and $L \to_\cT K$. If $\cS$ is an $H$-transfer system, then $\to_{\ind^G_H \cS}$ is the reflexive and transitive closure of the binary relation 
    \[
    \{L^g \to K^g \mid g\in G \text{, } L \leq K \leq H \text{ and }L \to_\cS K\}.
    \]
\end{proposition}

We can find similar characterisations of $\coind^G_H$ if $G$ is an abelian group.

\begin{proposition}
    \label{prop:coindCharacterisation}
    Let $H \leq G$ be abelian groups, and let $\cS$ be an $H$-transfer system. Then $\coind^G_H\cS$ is determined by the rule that $L \to_{\coind^G_H\cS} K$ if and only if $L\cap H \to_\cS K \cap H$.
\end{proposition}
\begin{proof}
    It is first necessary to verify that the binary relation $\to_{\coind^G_H\cS}$ defined by the above rule indeed defines a transfer system. Note that closure under conjugacy is guaranteed since $G$ is assumed to be abelian. The other requirements of a transfer system follow from the fact that $\cS$ is a transfer system.
    
    Now we will check that the map $\coind^G_H$, as defined above, is a right adjoint to $\res^G_H$. To this end, suppose that $\cT$ is a $G$ transfer system such that $\cT \leq \coind^G_H \cS$. Then, if $L \leq K \leq H$ and $L \to_\cT K$, it follows that $L\to_{\coind^G_H \cS} K$. Since $K \leq H$, we have that $L\cap H = L $ and $K \cap H = K$ thus, $L \to_\cS K$. Therefore, $\res^G_H\cT \leq \cS$. Conversely, if $\res^G_H\cT \leq \cS$, and $L \leq K$ are subgroups of $G$ such that $L \to_\cT K$, then we have $L\cap H \to_\cT K \cap H$ by closure under restriction. Thus, $L \cap H \to_\cS K \cap H$ and $L \to_{\coind^G_H\cS} K$, meaning $\cT \leq \coind^G_H \cS$ as required.
    \end{proof}

    \begin{remark}
        Our construction of $\coind^G_H$ becomes more complicated when one moves to non-abelian groups. This is because the rule that $L \to_{\coind^G_H \cS} K$ if and only if $L\cap H \to_\cS K \cap H$ does not produce a order relation that is closed under conjugacy in general. One would assume that we simply take the transfer system generated by this relation. However, we are constructing a right adjoint, and must instead take the largest transfer system contained in the relation we have described, see \cite[Proposition 5.9]{rubin2021operadic} for more details. Since the present paper only requires coinduced transfer systems for abelian groups, we have omitted the non-abelian case to avoid further complications.
    \end{remark}
    
    Suppose $N$ is a normal subgroup of $G$. If $\pi: G \to G/N$ is the quotient map, then $\pi_L:\trans(G) \to \trans(G/N)$ is known as the fixed-points map and is denoted by $(-)^N$. See \cite[Appendix B]{blumberg2015operadicmultiplicationsequivariantspectra} for a discussion of the corresponding fixed-point operads. The map $\pi_L^{-1}$ will be called inflation and denoted by $\inf^{G}_{G/N}$. The following characterisations can be deduced from \cite[Lemma 5.7]{rubin2021operadic}.

    \begin{proposition}
    \label{prop:fixedpointTransfer}
        If $\cT$ is a $G$ transfer system, then $\cT^N$ is determined by the rule that $K/N \to_{\cT^N} H/N$ if and only if $N \leq K \leq H$ and $K \to_{\cT} H$. If $\cS$ is a $G/N$-transfer system, then $\inf^G_{G/N}\cS$ is determined by the rule that $M \to_{\inf^G_{G/N}\cS} L$ if there exist $K,H \leq G$ such that $N \leq K \leq H$, $L \leq H$, $M = K \cap L$ and $K/N \to_\cS H/N$.
    \end{proposition} 
    In other words, the $N$ fixed point transfer system $\cT^N$ remembers only transfers in $\cT$ above $N$, and the inflation $\inf^G_{G/N}\cS$ is the closure of the relation $\to_\cS$ under restriction.
    
    We shall check directly which of these change-of-group maps preserve the saturated and disc-like properties, making use of the following lemma, whose proof we leave to the reader. 

    \begin{lemma}
     \label{lem:modularProperty}
        If $G$ is a group, $H$ is a subgroup of $G$, and $N$ is a normal subgroup of $G$, then $NH$ is the smallest subgroup of $G$ which contains both $N$ and $H$. Furthermore, the following two statements hold: 
        \begin{enumerate}
            \item If $K \geq N$, then $NH \cap K=N(H\cap K)$.
            \item If $K \geq H$, then $NH \cap K=(N \cap K)H$. \qed
        \end{enumerate} 
    \end{lemma}
    
    \begin{proposition}
        \label{prop:changeDiscSaturated}
        Let $G$ be a finite group, let $H\leq G$, $N \unlhd G$ and let $\cT \in \trans(G)$, $\cS \in \trans(H)$, $\cR \in \trans(G/N)$. Then, the following statements hold: 
        \begin{enumerate}
            \item If $\cT$ is saturated, then $\res^G_H\cT$ and $\cT^N$ are saturated.
            \item If $\cR$ is saturated, then $\inf^G_{G/N} \cR$ is saturated. 
            \item If $\cT$ is disc-like, then $\res^G_H\cT$ and $\cT^N$ are disc-like.
            \item If $\cR$ is disc-like, then $\inf^G_{G/N} \cR$ is disc-like.
            \item If $G$ is abelian and $\cS$ is disc-like, then $\coind^G_H \cS$ is disc-like. 
        \end{enumerate}
    \end{proposition}
    \begin{proof}
        (1) If $\cT$ is saturated, then any triangle $M \leq L \leq K$ of subgroup inclusions satisfies the two-out-of-three property. In particular, this is true if $K \leq H$ and $K \geq N$. It follows immediately that $\res^G_H\cT$ and $\cT^N$ are saturated. 

        (2) Recall that $\inf^G_{G/N} \cR$ is obtained by closing $\cR$ under restriction. Thus, if we suppose that $M \to_{\inf^G_{G/N} \cR} K$ and $M \leq L \leq K$, we wish to show that $L \to_{\inf^G_{G/N} \cR} K $. Note that $M \to_{\inf^G_{G/N} \cR} K$ implies there exist subgroups $A , C \leq G$ such that $N \leq A \leq C$, $A/N \to_\cR C/N$, $K \leq C$ and $A \cap K=M$. By Proposition \ref{prop:saturatedNormalCore}, we may assume that $A$ is normal in $C$. If we define $B = AL$, then we have $B/N \to_\cR C/N$ since $\cR$ is saturated. By Lemma \ref{lem:modularProperty}, $B \cap K= AL \cap K = (A \cap K) L = L$. Thus, $L \to_{\inf^G_{G/N} \cR} K$ by restriction of $B \to C$ to $K$. 
         
        (3) Let $\cT$ be disc-like. It is evident that $\cT^N$ is disc-like; it is generated by the relation $\{H/N \to G/N \mid H \to_\cT G\}$. To see that $\res^G_H\cT$ is disc-like, we use Proposition \ref{prop:disclikeRestriction}. In particular, if $\cT$ is disc-like and $L \to_{\res^G_H \cT} K$, then there exists a subgroup $H' \leq G $ such that $H' \to_\cT G$ and $H' \cap K = L$. By restriction, we have $H' \cap H \to_{\res^G_H\cT} H$ and $H' \cap H \cap K = H' \cap K = L$, so $L \to K$ is obtained by restriction from a transfer with target $H$.  

        (4) Let $\cR$ be disc-like, since $\inf^G_{G/N} \cR$ is obtained by closing $\cR$ under restriction, we see that $\inf^G_{G/N} \cR$ is generated by $\{H \to G \mid H/N \to_\cR G/N\}$, making $\inf^G_{G/N} \cR$ disc-like.

        (5) Now suppose $G$ is an abelian group and $\cS$ is disc-like as in the fourth statement. Recall that $L \to_{\coind^G_H\cS} K$ if and only if $L \cap H \to_\cS K \cap H$. Since $\cS$ is disc-like, there exists a relation $M \to_{\cS} H$ such that $M \cap (K \cap H)= M \cap K = L \cap H$. By Lemma \ref{lem:modularProperty}, $ML \cap H=M(L \cap H)=M$. We see that $ML \to_{\coind^G_H\cS} G$. Furthermore, $ML\cap K=L(M\cap K)=L$, thus $L \to K$ is obtained by restriction of $ML \to G$ to $K$, making $\coind^G_H\cS$ disc-like.
    \end{proof}
    
    For lattices of universes, there are similar change-of-group maps $\ind^G_H$ and $\res^G_H$, arising from the induction and restriction functors between the categories of $H$-representations and $G$-representations. Similarly we have maps $(-)^N$ and $\inf^G_{G/N}$ arising from the fixed points and inflation of universes. These are compatible with the change-of-group maps between lattices of transfer systems in the following ways. 

    \begin{proposition}
        \label{prop:changeofgroupUniverses}
        Let $\mathbb{K}$ be either $\R$ or $\C$, then:
        \begin{enumerate}
            \item The maps $\cT_\cD,\cT_\cL : \mathbf{Uni}_\mathbb{K}(G) \to \trans(G)$ commute with restriction.
            \item If $N$ is a normal subgroup of $G$ and $[\cU] \in \mathbf{Uni}_{\mathbb{K}}(G)$, then $(\cT_{\cD(\cU)})^N=\cT_{\cD(\cU^N)}$.
            \item If $[\cW] \in \mathbf{Uni}_{\mathbb{K}}(G/N)$, then $\cT_{\cD(\inf^G_{G/N}\cW)}=\inf^G_{G/N}\cT_{\cD(\cW)}$.
            \item If $G$ is an abelian group, $H$ is a subgroup of $G$ and $[\cU] \in \mathbf{Uni}_\mathbb{K}(H)$, then $\cT_{\cD(\ind^G_H \cU)}=\coind^G_H\cT_{\cD(\cU)}$.
        \end{enumerate}
    \end{proposition}
    \begin{proof}
        (1) The first part of the proposition immediately follows from the definitions of the transfer systems $\cT_{\cD(\cU)}$ and $\cT_{\cL(\cU)}$, and the fact that the restrictions of representations are functorial with respect to subgroup inclusions. 

        (2) By Proposition \ref{prop:fixedpointTransfer}, the second part of the proposition is equivalent to stating that $\frac{(H/N)}{(K/N)}$ embeds $H/N$-equivariantly into $\cU^N$ if and only $H/K$ embeds $H$-equivariantly into $\cU$. This holds by the canonical bijection $H/K \cong \frac{(H/N)}{(K/N)}$.

        (3) We can deduce from Proposition \ref{prop:fixedpointTransfer} that $H \to_{\inf^G_{G/N}\cD(\cW)} G$ if and only if $N\leq H$ and $H/N \to_{\cD(\cW)}G/N$. Thus, $H \to_{\inf^G_{G/N}\cD(\cU)} G$ if and only if $\frac{(G/N)}{(H/N)}$ embeds $G/N$-equivariantly into $\cW$. This occurs if and only if $G/H$ embeds $G$-equivariantly into $\inf^G_{G/N}\cW$. By Proposition \ref{prop:changeDiscSaturated}, both of the transfer systems $\cT_{\cD(\inf^G_{G/N} \cW)}$ and $\inf^G_{G/N}\cT_{\cD(\cW)}$ are disc-like. Thus, we can conclude the desired result, that $\cT_{\cD(\inf^G_{G/N} \cW)}=\inf^G_{G/N}\cT_{\cD(\cW)}$.

        (4) Recall that $\ind^G_H \cU$ can be viewed as a $G/H$-indexed direct sum of $\cU$. Hence, a point in $\ind^G_H \cU$ is a formal sum $\sum_{[g] \in G/H} gx_g$, where each $x_g$ is a point in $\cU$. Since $G$ is abelian, the $H$ action on such a sum is diagonal, so $\res^G_H \ind^G_H \cU \cong \cU^{\oplus [G:H]} \cong \cU$ and it follows that $\res^G_H(\cT_{\cD(\ind^G_H \cU)})=\cT_{\cD(\cU)}$. Since $\coind^G_H$ is a right adjoint to $\res^G_H$ on transfer systems, then $\cT_{\cD(\ind^G_H \cU)} \leq \coind^G_H \cT_{\cD(\cU)}$. To see the other inclusion, note that both $\coind^G_H \cT_{\cD(\cU)}$ and $\cT_{\cD(\ind^G_H \cU)}$ are disc-like, thus it suffices to see that $K \to_{\coind^G_H \cD(\cU)} G$ implies that $K \to_{\cD(\ind^G_H \cU)} G$. Recall that $K \to_{\coind^G_H \cD(\cU)} G$ implies $K \cap H \to_{\cD(\cU)} H$, which in turn means that there exists $x \in \cU$ with stabiliser $H_x = H \cap K$. We wish now to construct a point in $\ind^G_H \cU$ with stabiliser $K$. Such a point is given by the sum $\sum_{[k] \in KH/H} kx$.
    \end{proof}
    
    \begin{remark}
        To address the possible confusion about the fourth part of the above proposition, we note that induced and coinduced representations of finite groups are naturally isomorphic, hence the resulting maps on lattices of universes coincide and we have only written $\ind^G_H$. This sort of ambidexterity is not true on the level of transfer systems, hence the mismatch.
    \end{remark}

\section{Equivalences between linear isometries and little discs operads.}
\label{section:equivalence}
It has been known for some time that equivariant linear isometries and little discs operads need not be equivalent in general, even over a fixed universe $\cU$. However, it was not known until now for which universes $\cU$ there exists a zig-zag of weak equivalences connecting $\cL(\cU)$ and $\cD(\cU)$. In \cite[\S 4]{blumberg2015operadicmultiplicationsequivariantspectra}, Blumberg and Hill showed that such a zig-zag exists when $\cU$ is isomorphic to $\R[G/N]^\infty$ for some normal subgroup $N$; the infinite direct sum of the regular representation of $G/N$. In this section, we use the language of transfer systems to prove that this is the only case where $\cT_{\cL(\cU)}=\cT_{\cD(\cU)}$. We will begin by examining what information about a universe $\cU$ can be recovered from the minimal fibrant subgroups of $\cT_{\cD(\cU)}$ and $\cT_{\cL(\cU)}$.

\begin{lemma}
    \label{lemma:linearMinfibrant}
    Let $\cU$ be a $G$-universe, then: 
    \begin{itemize}
        \item If $N$ is the minimal fibrant subgroup for $\cT_{\cL(\cU)}$, then $\cU^N \cong \R[G/N]^\infty$.
        \item If $N$ is the minimal fibrant subgroup for $\cT_{\cD(\cU)}$, then $N$ acts trivially on $\cU$.
    \end{itemize}
\end{lemma}
\begin{proof}
    Suppose that $N$ is the minimal fibrant subgroup for $\cT_{\cL(\cU)}$. Then $N \to_{\cL(\cU)} G$, hence there is a $G$ equivariant embedding $$ \ind^G_N\res^G_N \cU \to \cU $$
    since $\cU$ contains a trivial summand, then the universe $\ind^G_N\res^G_N \cU$ contains the subrepresentation $\ind^G_N\res^G_N \R \cong \R[G/N]$. Thus, $\cU$ contains $\R[G/N]$, making $\cU^N$ a complete universe for $G/N$, meaning $\cU^N \cong \R[G/N]^\infty$.

    Suppose now that $N$ is the minimal fibrant subgroup for $\cT_{\cD(\cU)}$. Then we have a $G$-equivariant embedding $G/N \to \cU$, but for all $H < N$ there is no $G$-equivariant embedding $G/H \to \cU$. That is, whenever $H < N$, the space of $H$-fixed points in $\cU$ is equal to the space of $N$-fixed points in $\cU$. By taking $H=\{e\}$, we see that $N$ fixes all points in $\cU$.
\end{proof}

The above lemma turns out to be sufficient information to prove our desired result. 

\begin{theorem}
    \label{thm:discEquivLinear}
    Let $G$ be a finite group and $\cU$ be a $G$-universe. Then $\cL(\cU) \simeq \cD(\cU)$ if and only if there exists a normal subgroup $N \unlhd G$ such that $\cU \cong \R[G/N]^\infty$.
\end{theorem}
\begin{proof}
    The `if' part of the statement is \cite[Theorem 4.21]{blumberg2015operadicmultiplicationsequivariantspectra}. Here we give a purely combinatorial proof using the following sequence of refinements of transfer systems: 
    \[ 
    \cT_{\cL(\R[G/N]^\infty)} \leq  \cT_{\cD(\R[G/N]^\infty)} \leq \langle H \to G\mid N \leq H \rangle \leq \cT_{\cL(\R[G/N]^\infty)},
    \]
    where the first refinement follows from Theorem \ref{thm:discsContainsLinear}. For the second, note that since $\cT_{\cD(\R[G/N]^\infty)}$ is disc-like, it suffices to show that whenever $G/H$ embeds into $\R[G/N]^\infty$, then $H$ contains $N$. Let $K$ be a subgroup which does not contain $N$ and suppose, for the sake of contradiction, that $K \to_{\cD(\R[G/N]^\infty)} G $. By restriction to $N$ it follows that $K \cap N\to_{\cD(\R[G/N]^\infty)} N $. However, the subgroup $K \cap N$ is strictly contained in $N$ and $N$ acts trivially on $\R[G/N]^\infty$. Thus, there can be no $N$-equivariant embedding of $N/(K\cap N)$ into $\R[G/N]^\infty$, giving the required contradiction. For the final refinement, first note that $$\ind^G_N \res^G_N \R[G/N]=\ind^G_N \R^{[G:N]}=\R[G/N]^{[G:N]}.$$ As a result, $N \to_{\cL(\R[G/N]^\infty)}G$. Since $\cT_{\cL(\R[G/N]^\infty)}$ is saturated, we see that there is a transfer $H \to_{\cL(\R[G/N]^\infty)} G$ for all $H \geq N$. In particular, $\cT_{\cL(\R[G/N]^\infty)}$ contains the transfer system generated by the relation $\{H \to G \mid N\leq H \}$. 

    We will now prove `only if'. Observe that since $\cL(\cU) \simeq\cD(\cU)$, the transfer systems $\cT_{\cL(\cU)}$ and $\cT_{\cD(\cU)}$ have the same minimal fibrant subgroup $N$. By Lemma \ref{lemma:linearMinfibrant}, $N$ acts trivially on $\cU$ and $\cU^N \cong \R[G/N]^\infty$. Hence, $\cU = \cU^N \cong \R[G/N]^\infty$ as desired.
\end{proof}

    Theorem \ref{thm:discEquivLinear} is not the end of the story, however. We have not yet addressed the case in which $\cL(\cU)$ and $\cD(\cW)$ are equivalent for distinct universes $\cU$ and $\cW$. Indeed, this is still possible. For example, if $G = C_5$ and $V$ is any irreducible $C_5$ representation, then $\cL((\R \oplus V)^\infty)$ and $\cD(\R^\infty)$ both give the trivial transfer system (containing only reflexive order relations), but $(\R \oplus V)^\infty$ is not isomorphic to $\R^\infty$ since the former has a non-trivial action by $G$. Note that if $\cT =\cT_{\cL(\cU)} = \cT_{\cD(\cW)}$, then $\cT$ must be both saturated and disc-like. 
    
    It turns out that transfer systems $\cT$ which are both saturated and disc-like are rather simple to classify.
    
    \begin{definition}
        \label{def:bisaturated}
         A transfer system is \textbf{bisaturated} if it is both saturated and disc-like. We denote the poset of bisaturated $G$ transfer systems by $\mathrm{BiSat}(G)$.
    \end{definition}
    
    \begin{theorem}
        \label{thm:bisaturated}
        The map taking a transfer system $\cT$ to its minimal fibrant subgroup restricts to an order-reversing bijection
        \[
           \mathrm{BiSat}(G) \to \sub_{\unlhd}(G).
        \]
    \end{theorem}
    \begin{proof}
        Surjectivity follows immediately from the fact that for any normal subgroup $N$, the bisaturated transfer system $\cT_{\cD(\R[G/N]^\infty)}=\cT_{\cL(\R[G/N]^\infty)}$ has minimal fibrant subgroup $N$.

        For injectivity, suppose that $\cT$ and $\cS$ are bisaturated transfer systems which have the same minimal fibrant subgroup $M$. We certainly have that $M \to_\cP G$. Since $\cS$ is saturated, it follows that $H \to_\cP G$ for all subgroups $H \geq M$. We also have that $\cT$ is disc-like, hence it is generated by some set of transfers $\{H_i \to G\}_{i\in I}$ where each $H_i$ necessarily contains $M$. Since $\cS$ possesses all such transfers, it follows that $\cT \leq \cS$. By symmetry, we can conclude that $\cT=\cS$.
    \end{proof}

    \begin{remark}
        The notion of a bisaturated transfer system has recently been formulated independently in \cite{klanderman2025characterizingtransfersystemsnonabelian}, where a conjectural connection is made between bisaturated transfer systems and minimal chains in the lattice of transfer systems. We hope that the results of this section will be useful for this conjecture.
    \end{remark}

    Applying Lemma \ref{lemma:linearMinfibrant} and Theorem \ref{thm:bisaturated} gives the following operadic corollary.

    \begin{corollary}
        If a transfer system $\cT$ is bisaturated, then there exists a normal subgroup $N$ such that $\cT = \cT_{\cD(\R[G/N]^\infty)}=\cT_{\cL(\R[G/N]^\infty)}$. In particular, if $\cU$ and $\cW$ are $G$ universes, and there exists a zig-zag of weak equivalences connecting $\cL(\cU)$ and $\cD(\cW)$, then there exists a normal subgroup $N$ and a sequence of embeddings $\cW \hookrightarrow \R[G/N]^\infty \hookrightarrow \cU$ such that $\cL(\cU) \simeq \cL(\R[G/N]^\infty)\simeq \cD(\R[G/N]^\infty)\simeq \cD(\cW)$. \qed
    \end{corollary}

\section{The maximality of linear isometries operads over abelian groups}
\label{section:maximality}

Since the map $\cL(\cU) \to \cD(\cU)$, in the homotopy category of $N_\infty$-operads, is not an isomorphism in general, it is natural to wonder how else we can characterise the relationship between $\cL(\cU)$ and $\cD(\cU)$. It is known that we can replace $\cD(\cU)$ with the equivariant infinite Steiner operad $\cK(\cU)$ up to an explicit weak equivalence, and that the pair $(\cL(\cU),\cK(\cU))$ forms what is known as an operad pair, see \cite{guillou2017equivariant}.

Operad pairs were first introduced in \cite[\S 9]{may1977infinite} in the non-equivariant setting. The definition naturally generalises to the equivariant setting as follows. 

\begin{definition}
\label{def:operadPair}
        Two $G$-operads $\cO$ and $\cP$ form an \textbf{operad pair} $(\cP , \cO)$ if there exist $G$-equivariant maps: 
        \[
        \cP(k) \times \cO(n_1) \times \cdots\times \cO(n_k) \to \cO(n_1 \times \cdots \times n_k)
        \]
        which are suitably compatible with the operad structures of $\cO$ and $\cP$.
\end{definition}

There is also the more combinatorial notion of a compatible pair of transfer systems.

\begin{figure}[ht!]
    \centering
\begin{tikzcd}
	& H &&&& H \\
	K && L & \Longrightarrow & K && L \\
	& {K \cap L} &&&& {K\cap L}
	\arrow[from=2-1, to=1-2]
	\arrow[from=2-5, to=1-6]
	\arrow[dashed, from=2-7, to=1-6]
	\arrow[dashed, from=3-2, to=2-1]
	\arrow[from=3-2, to=2-3]
	\arrow[dashed, from=3-6, to=2-5]
	\arrow[from=3-6, to=2-7]
\end{tikzcd}
    \caption{A pair of diagrams depicting the condition which a compatible pair of transfer systems $(\cS,\cT)$ must satisfy, where the solid arrows depict relations $\to_\cS$ and the dashed arrows relations $\to_\cT$.}
    \label{fig:compatiblePair}
\end{figure}
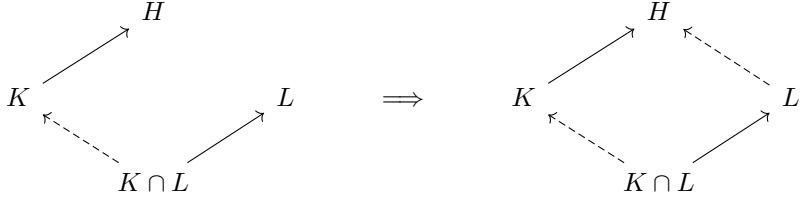

\begin{definition}[{\cite[Definition 4.6]{chan2024bi}}]
\label{def:compatiblePair}
    A pair of transfer systems $(\cS,\cT)$ is said to be \textbf{compatible} if the following condition holds: If $H$ is a subgroup of $G$ and $L, K \leq H$ are subgroups such that $K \to_\cS H$ and $K \cap L \to_\cT K$, then $L \to_\cT H$, see Figure \ref{fig:compatiblePair}.     
\end{definition}

It has been shown in \cite[Theorem 5.1]{chan2025realizingcompatiblepairstransfer} that if $(\cP, \cO)$ form an operad pair of $N_\infty$-operads, then the corresponding pair of transfer systems $(\cT_\cS,\cT_\cO)$ is compatible. In particular, the pair $(\cT_{\cL(\cU)},\cT_{\cD(\cU)})=(\cT_{\cL(\cU)},\cT_{\cK(\cU)})$ always forms a compatible pair. In \cite[Corollary 6.19]{blumberg2021}, Blumberg and Hill showed (in the equivalent context of indexing categories) that for each transfer system $\cT$, there exists a unique maximal transfer system $\Hull(\cT)$ with the property that $(\Hull(\cT),\cT)$ form a compatible pair. The pair $(\Hull(\cT),\cT)$ is accordingly said to be maximally compatible, and the transfer system $\Hull(\cT)$ is known as the multiplicative hull of $\cT$. 

\subsection{A closure operator on the lattice of universes}
\label{subsec:closureOperator}
In this section, we will show that for all abelian groups $G$ and all $G$ universes $\cU$, there exists a universe $\overline{\cU}$ such that $\cT_{\cL(\overline{\cU})}=\Hull(\cT_{\cD(\cU)})$.

\begin{definition}
    \label{def:universeClosure}
    Let $G$ be a finite group and let $\cU$ be a $G$-universe. Define the \textbf{(additive) closure} $\overline\cU$ of $\cU$ as the unique maximal universe in $\mathbf{Uni}(G)$ such that $\cT_{\cD(\cU)} = \cT_{\cD(\overline\cU)}$. Such a universe exists since $\cT_\cD: \mathbf{Uni}(G) \to \mathrm{Tr(G)}$ preserves meets. We say that a universe $\cU$ is \textbf{(additively) closed} if $\cU \cong \overline\cU$.
\end{definition}

\begin{lemma}
    \label{lem:closedProperties}
    Let $\cU$ be a $G$-universe, then the following statements are equivalent: 
    \begin{enumerate}
        \item $\cU$ is additively closed.
        \item For all $G$-universes $\cW$, if $\cT_{\cD(\cW)} \leq \cT_{\cD(\cU)}$ then $\cW$ embeds into $\cU$. 
        \item If $V$ is an irreducible $G$-representation, such that whenever $G/H$ embeds equivariantly into $V$, we have $H \to_{\cD(\cU)} G$, then $V$ embeds into $\cU$.
    \end{enumerate}
\end{lemma}
\begin{proof}
    We begin with the implication $(2) \Rightarrow (1)$. Note that, by the definition of the additive closure, we have an embedding $\cU \hookrightarrow \overline{\cU}$, and that $\cT_{\cD(\cU)} = \cT_{\cD(\overline{\cU})}$. By $(2)$, it follows that there is an embedding $\overline{\cU} \hookrightarrow \cU$. Hence, $\cU \cong \overline \cU$ as desired. 

    Now we prove that $(3)$ implies $(2)$. Recall that we can decompose a universe $\cW$ as $\cW \cong (\R \oplus \bigoplus_{i=0}^n V_i)^\infty$, where $\{V_1, \dots ,V_n\}$ are the irreducible subrepresentations of $\cW$. Since $\cT_{\cD(\cW)}\leq \cT_{\cD(\cU)}$, it follows from $(3)$ that each $V_i$ embeds into $\cU$. Hence we have an embedding $\cW \hookrightarrow \cU$.

    Finally, we prove that $(1)$ implies $(3)$. Let $\{V_1,\dots,V_n\}$ be a collection (up to isomorphism) of all irreducible $G$-representations such that whenever $G/H$ embeds into $V$, we have $H\to_{\cD(\cU)} G$. We claim that $(V_1 \oplus \cdots \oplus V_n)^\infty$ is a model for $\overline{\cU}$. Indeed, if $\cW$ is a $G$-universe with $\cT_{\cD(\cW)}=\cT_{\cD(\cU)}$, then each irreducible subrepresentation $W \subset \cW$ is contained in the list $\{V_1,\dots,V_n\}$ by assumption, hence we have an embedding $\cW \hookrightarrow (V_1 \oplus \cdots \oplus V_n)^\infty$. Thus, assuming $(1)$, we conclude that $\cU \cong \overline{\cU}\cong (V_1 \oplus \cdots \oplus V_n)^\infty$, so each $V_i$ embeds into $\cU$.
\end{proof}

\begin{remark}
    We have used this notation since the assignment of $\overline{\cU}$ to $\cU$ descends to a closure operator on the lattice $\mathbf{Uni}(G)$. This means that the following statements hold, where $[\cU] \in \mathbf{Uni}(G)$ is the isomorphism class of a universe $\cU$.
    \begin{itemize}
        \item $[\cU] \leq [\overline{\cU}]$.
        \item $[\cU] \leq [\cW]$ implies that $[\overline{\cU}] \leq [\overline{\cW}]$.
        \item $[\overline{\overline{\cU}}]=[\overline{\cU}]$.
    \end{itemize}
\end{remark}

\begin{remark}
    To justify our above use of the words additive and multiplicative, we note that compatible pairs were introduced by Blumberg-Hill \cite{blumberg2021} in the context of bi-incomplete Tambara functors, where one transfer system controls which multiplicative norms can exist and the other controls which additive transfers can exist. In particular, the transfer system $\cT_{\cD(\cU)}$ specifies the additive transfers in the homotopy groups of $G$-spectra indexed over $\cU$. We use the term `additive closure' because $\overline{\cU}$ is the largest universe which parameterises the same additive transfers as $\cU$.
\end{remark}

We shall now investigate the relationship between the closure operator $\overline{(-)}$ and various change-of-group maps.

\begin{proposition}
    If $G$ is a finite group, $\cU$ is a closed universe, and $N$ is a normal subgroup of $G$, then the $G$-universe $\cU^N$ is closed.
\end{proposition}
\begin{proof}
    Recall from Proposition \ref{prop:changeofgroupUniverses}, that $\cT_{\cD(\cU^N)}=(\cT_{\cD(\cU)})^N$ as transfer systems. Let $\cU$ be a closed $G$-universe and suppose that $\cW$ is a $G/N$ universe such that $\cT_{\cD(\cW)}=\cT_{\cD(\cU^N)}=(\cT_{\cD(\cU)})^N$. Note that we may view $\cW$ as a $G$-universe by inflation. Again, by Proposition \ref{prop:changeofgroupUniverses}, $\cT_{\cD(\inf^G_{G/N} \cW)} = \inf^G_{G/N}\cT_{\cD(\cW)}$, and we have a refinement of transfer systems $\cT_{\cD( \inf^G_{G/N}\cW)}=\inf^G_{G/N}\cT_{\cD(\cW)} \leq \cT_{\cD(\cU)}$ by the adjunction $\inf^G_{G/N} \dashv (-)^N$ on transfer systems. Since $\cU$ is closed, it follows from Lemma \ref{lem:closedProperties} that there exists an embedding $\inf^G_{G/N}\cW \hookrightarrow \cU$. By taking fixed points we get an embedding $\cW\cong (\inf^G_{G/N}\cW)^N \hookrightarrow \cU^N$, meaning that $\cU^N$ is closed.
\end{proof}

In what follows, we will work with abelian groups, since their representation theory and subgroup structure are particularly easy. As in \cite{rubin2021detecting}, a subgroup $H \leq G$ is cocyclic if the quotient $G/H$ is cyclic. For an abelian group $G$, the irreducible $G$-representations are precisely the irreducible representations over the quotients $G/H$, where $H$ is cocyclic, see \cite[\S 4.3]{rubin2021detecting}. The proposition below is an immediate consequence of this observation. 
\begin{proposition}
    \label{prop:closedCocyclic}
    If $G$ is an abelian group, and $\cU$ is a $G$-universe, then $\cU$ is closed if and only if for each cocyclic subgroup $H$, the $G/H$-universe $\cU^H$ is closed. 
\end{proposition}
\begin{proof}
    The `if' part of the proposition is just a special case of the previous proposition. For the `only if', we use Lemma \ref{lem:closedProperties}. Suppose that $V$ is an irreducible $G$-representation such that whenever $G/K$ embeds equivariantly into $V$, we have $K \to_{\cD(\cU)} G$. By our previous discussion, $V$ is the inflation of an irreducible representation of some cyclic quotient $G/H$. Hence whenever $G/K$ embeds equivariantly into $V$, then $H \leq K$, and $K \to_{\cD(\cU^H)} G$. Since $\cU^H$ is closed it follows from Lemma \ref{lem:closedProperties} that there is an embedding $V \hookrightarrow \cU^H \hookrightarrow \cU$. Again, by Lemma \ref{lem:closedProperties} we are done.  
\end{proof}

Let $C_n$ denote the cyclic group of order $n$. As noted by Rubin, each irreducible $C_n$-representation embeds into a $2$-dimensional representation $\lambda_n(m)$, where the generator $g$ acts as $e^{\frac{2\pi i m}{n}}$ on $\C$, thought of here as $\R^2$. By \cite[Lemma 5.13]{rubin2021detecting}, these have the following properties. 
    \begin{itemize}
        \item If $m\equiv m' \;\mathrm{mod} \;n $, then $\lambda_n(m)=\lambda_n(m')$.
        \item There are isomorphisms $\lambda_n(m)=\lambda_n(-m)$.
        \item If $d \mid n$, then $\res^{C_n}_{C_d} \lambda_n(m)=\lambda_d(m)$.
    \end{itemize}
As a result of this, every $C_n$-universe $\cU$ can be uniquely written in the form $\cU \cong (\bigoplus_{i \in I} \lambda_n(i))^\infty$ where $I \subseteq C_n$ is a subset which contains zero and is closed under additive inverses, see \cite[Lemma 5.14]{rubin2021detecting}. Following \cite{rubin2021detecting} and \cite{hafeez2022saturatedlinearisometrictransfer}, such sets are known as index sets and we denote this universe by $\cU_I$. It follows that $\res^{C_n}_{C_d} \cU_I \cong \cU_{(I \, \mathrm{mod} \, d)}$.

\begin{proposition}
    For an index set $I \subseteq \Z/n\Z$, the universe $\cU_I$ is closed if and only if whenever $C_d \to_{\cD(\cU_I)} C_n$, and $i \leq n$ is a natural number such that $\gcd(i,n)=d$, then $i \in I$. An index set is said to be closed if it satisfies this property.p
\end{proposition}
\begin{proof}
     By Lemma \ref{lem:closedProperties}, and since the irreducible representations of $C_n$ are free away from the origin, $\cU_I$ is a closed $C_n$ universe if and only if for each transfer of the form $C_d \to_{\cD(\cU_I)} C_n$, the universe $\cU_I$ contains every irreducible $C_n$-representation $V$ where $C_d$ appears as a stabiliser of a point in $V$. We shall now examine the stabilisers of points in $\lambda_n(m)$. Note that $e^{\frac{2\pi i m}{n}}$ acts non-trivially on all non-zero points in $\R^2$ as long as $m/n$ is not an integer. It follows that the stabiliser of all non-zero points in $\lambda_n(m)$ is the subgroup $C_{n/a} \leq C_n$ where $a$ is the smallest divisor of $n$ such that $ma/n$ is an integer. Specifically, we have $a=\frac{n}{\mathrm{gcd(m,n)}}$. Thus, there is a transfer $C_d \to_{\cD(\cU_I)} C_n$ if and only if $I$ contains every $i$ with the property that $\gcd(i,n)=d$, as desired.
\end{proof}

By the above proposition, we may reduce our study of closed universes over cyclic groups to elementary number theory. 

\begin{proposition}
    \label{prop:numberTheory}
    Let $n$ be a natural number and let $a$ and $d$ be divisors of $n$ such that $a$ and $d$ are coprime. If $i < d$ is coprime to $d$, then there exists $m \leq n$ such that $m \equiv i \,\mathrm{mod}\, d$ and $\gcd(m,n)=a$.
\end{proposition}
\begin{proof}
    Let $k$ be the natural number such that $n=kd$, and let $j=\frac{k}{\gcd(k,d)}$. Note that $j$ is the largest divisor of $n$ that is coprime with $d$. By Sunzi's theorem (usually known as the Chinese remainder theorem), there exists a natural number $m$ such that $m \equiv i \, \mathrm{mod} \, d$ and $m \equiv a \, \mathrm{mod} \, j$. It follows that $m$ and $d$ are coprime, which implies that $\gcd(m,n)=\gcd(m,j)=a$ as required.
\end{proof}

\begin{lemma}
    \label{lem:resClosed}
    Let $G$ be an abelian group and let $H$ be a subgroup of $G$. If $\cU$ is a closed $G$-universe, then $\res^G_H \cU$ is a closed $H$-universe.
\end{lemma}
\begin{proof}
    We first show this in the case that $G \cong C_n$ and $H \cong C_d$ for some divisor $d$ of $n$. Here, $\cU$ is isomorphic to $\cU_I$ for some index set $I$, by \cite[Lemma 5.14]{rubin2021detecting}. Since $\cU$ is assumed to be closed, then $I$ must be a closed index set. The statement of the lemma is then simply that $(I \, \mathrm{mod} \, d)$ is a closed index set. To prove this, suppose that $C_k \to_{\res^{C_n}_{C_d}\cU_I} C_d$, and let $i<d$ be a natural number such that $\gcd(i,d)=k$. We wish to show that $i \in I \, \mathrm{mod}\, d$. By Proposition \ref{prop:disclikeRestriction}, the transfer $C_k \to_{\res^{C_n}_{C_d}\cU_I} C_d$ can be obtained by restriction of a transfer $C_a \to_{\cD(\cU_I)} C_n$ to $C_d$. It follows that $C_k = C_d \cap C_a = C_{\gcd(d,a)}$. We may now divide each present number by $k$ to assume that $i$ and $d$ are coprime and that $d$ and $a$ are coprime. Applying Proposition \ref{prop:numberTheory} gives a natural number $m$ such that $m \equiv i \, \mathrm{mod} \, d$ and $\gcd(m,n)=a$. If needed, we may now multiply everything by $k$ to return each number to its original value. Note that since $C_a \to_{\cD(\cU_I)} C_n$ and $I$ is a closed index set, then $m \in I$. We can conclude that $i \in I \, \mathrm{mod} \, d$. 

    We will now show that the general abelian case of the present lemma reduces to the cyclic case we have just discussed. Let $H\leq G$ be a pair of abelian groups and assume that $\cU$ is a closed $G$-universe. By Proposition \ref{prop:closedCocyclic}, it suffices to show that $(\res^G_H \cU)^K$ is closed for every cocyclic subgroup $K \leq H$. Consider the abelian group $G/K$ and note that $(\res^G_H \cU)^K \cong \res^{G/K}_{H/K}(\cU^K)$. If $G/K$ is cyclic, then we have already proven that $\res^{G/K}_{H/K}(\cU^K)$ is closed. Otherwise, we may use the structure theorem for finite abelian groups to find a cocyclic subgroup $K' \leq G$ such that $K' \cap H=K$. The composite $H \to G \to G/K'$ then induces an injective group homomorphism $ \phi: H/K \to G/K'$. If $H'$ is the subgroup of $G$ such that $H'/K'$ is the image of $\phi$, then we have a commutative diagram as below. 
\[\begin{tikzcd}[cramped]
	G && {G/K} \\
	{H'} & {H'/K'} & {H/K}
	\arrow[two heads, from=1-1, to=1-3]
	\arrow[hook, from=2-1, to=1-1]
	\arrow[two heads, from=2-1, to=2-2]
	\arrow["\cong", from=2-2, to=2-3]
	\arrow[hook, from=2-3, to=1-3]
\end{tikzcd}\]
     It can then be shown that $(\res^G_{H'} \cU)^{K'}$ corresponds to $\res^{G/K}_{H/K} (\cU^K)$ under the isomorphism $H'/K' \cong H/K$, see \cite[Lemmas 2.2 \& 2.5]{li2025non}. Note that $(\res^G_{H'} \cU)^{K'}\cong \res^{G/K'}_{H'/K'}(\cU^{K'})$. Since $\cU^{K'}$ is closed and $H'/K' \to G/K'$ is an inclusion of cyclic groups, the $H'/K'$ universe $\res^{G/K'}_{H'/K'}(\cU^{K'})$ is closed. By passing through the isomorphism $H'/K' \cong H/K$, we see that $(\res^G_H \cU)^K\cong \res^{G/K}_{H/K} (\cU^K)$ is closed, completing the proof.
\end{proof}

\subsection{Computing the multiplicative hull of infinite little discs operads}

We now turn to a discussion of maximal compatible pairs of transfer systems. If $G$ is an abelian group and $\cT$ is a $G$ transfer system, then the following definition gives a new characterisation of $\Hull(\cT)$.

\begin{definition}
    Given an abelian group $G$ and a transfer system $\cT$, define a binary relation $\to_{\cH(\cT)}$ on $\sub(G)$ by the rule that $K \to_{\cH(\cT)} H $ if and only if $K \leq H \leq G$ and $\coind^H_K \res^G_K \cT\leq \res^G_H\cT$, where $\coind^H_K$ and $\res^G_K$, and $\res^G_H$ are the coinduction and restriction maps discussed in Section \ref{subsec:changeGroup}.
\end{definition}

\begin{theorem}
    The binary relation $\to_{\cH(\cT)}$ defines a transfer system $\cH(\cT)$.
\end{theorem}
\begin{proof}
    Let us check the requirements of a transfer system. Refinement of inclusions holds by definition, closure under conjugacy is trivial since $G$ is assumed to be abelian, and reflexivity is immediate since $\coind^H_H\res^G_H \cT= \res^G_H \cT$. Now, for transitivity, suppose that $L \to_{\cH(\cT)} K$ and $K \to_{\cH(\cT)}H$. Then, the inequalities $\coind^K_L \res^G_L \cT \leq \res^G_K \cT$ and $\coind^K_L \res^G_L \cT \leq \res^G_H \cT$ give rise to the following sequence of refinements, showing that $L \to_{\cH(\cT)} H$: 
    \[
    \coind^H_L \res^G_L \cT = \coind^H_K \coind^K_L \res^G_L \cT \leq \coind^H_K \res^G_K \cT \leq \res^G_H \cT.
    \]
    It remains to check closure under restriction. Suppose that $L \to_{\cH(\cT)} H$ and $K \leq H$. By Proposition \ref{prop:coindCharacterisation}, we see that the transfer system $\coind^K_{K\cap L}\res^G_{K\cap L} \cT$ consists of subgroup inclusions $A \to B$ such that $A \leq B \leq K$ and $A \cap (K \cap L) \to_{\res^G_{K\cap L} \cT} B \cap(K \cap L)$. Note that $A \cap K \cap L= A\cap L$ and $B \cap K \cap L= B\cap L$, hence $A\cap L \to_{\res^G_L \cT} B \cap L$. Thus, $A\to_{\res^H_K \coind^H_L \res^G_L \cT} B$, and we have the following sequence of refinements: 
    
    \[
    \coind^K_{K \cap L}\res^G_{K \cap L} \cT \leq \res^H_K \coind^H_L \res^G_L \cT \leq \res^H_K \res^G_H \cT = \res^G_K \cT,
    \]
    which shows that $L \cap K \to_{\cH(\cT)} K$, completing the proof.
\end{proof}

One should intuitively think of a transfer $K\to_{\cH(\cT)}H$ as encoding the existence of twisted morphisms $\prod_{hK \in H/K}\cT \to \cT$ in the homotopy category of operads (using `twisted' to mean $\Gamma$-equivariant for a graph subgroup $\Gamma \leq G \times \Sigma_{[H:K]}$ encoding the $H$ action on $H/K$). However, it is currently unclear how to make this intuition precise.

\begin{theorem}
\label{thm:internalLinearisHull}
    Let $G$ be an abelian group and let $\cT$ be a $G$ transfer system. Then, $\Hull(\cT)=\cH(\cT)$.
\end{theorem}
\begin{proof}
    We shall first prove that the pair $(\cH(\cT),\cT)$ is a compatible pair of transfer systems, then we will prove that it is maximal.

    Let us check that $(\cH(\cT),\cT)$ is a compatible pair of transfer systems. First we must check that $\cH(\cT) \leq \cT$. For this, let $K \leq H \leq G$ be subgroups such that $K \to_{\cH(\cT)} H$. Then, $\coind^H_K\res^G_K \cT \leq \res^G_H \cT$. Note that $K \to_{\coind^H_K\res^G_K \cT} H $ must hold by Proposition \ref{prop:coindCharacterisation}, thus we also have $K \to_{\res^G_H\cT} H$. By Proposition \ref{prop:resindCharacterisation}, we have $K \to_\cT H$ as required. For the second condition of Definition \ref{def:compatiblePair}, suppose that $L,K \leq H \leq G$ are subgroups of $G$ such that $K \to_{\cH(\cT)} H$ and $K \cap L \to_\cT H$. By Proposition \ref{prop:coindCharacterisation}, $L \to_{\coind^H_K \res^G_K \cT} H$, since $L \cap K \to_{\res^G_K \cT } K$. Furthermore, $\coind^H_K \res^G_K \cT \leq \res^G_H \cT$, hence $L \to_\cT H$.

\[\begin{tikzcd}[cramped]
	L & H \\
	{L\cap K} & K
	\arrow["\cT"', dashed, from=1-1, to=1-2]
	\arrow["{\cH(\cT)}", from=2-1, to=1-1]
	\arrow["\cT", dashed, from=2-1, to=2-2]
	\arrow["{\cH(\cT)}"', from=2-2, to=1-2]
\end{tikzcd}\]

    Finally, we will show that the compatible pair $(\cH(\cT),\cT)$ is maximal. If $L\leq M$ and $L\to_\cT M$ is a transfer which is not in $\cH(\cT)$, then $\coind^M_L\res^G_L \cT \nleq \res^G_M \cT$. Let $\cS$ be a transfer system such that $L\to_\cS M$ and $\cS \leq \cT$. We wish to show that $(\cS,\cT)$ is not a compatible pair of transfer systems. Since $\coind^M_L\res^G_L \cT \nleq \res^G_M \cT$, there exist subgroups $K$ and $H$ such that $K \leq H \leq M$ and $K \to_{\coind^M_L\res^G_L\cT} H$ but $K \to H$ is not in $\res^G_M \cT$. Thus $K \cap L \to_\cT H\cap L$ but $K \to H$ is not in $\cT$. We also have $H \cap L \to_\cS H$ by restriction of $L \to_\cS M$ to $H$, and $K \cap L \to_\cS K$ by restriction of $L \to_\cS M$ to $K$. These subgroups can be arranged into a square as below:
    
\[\begin{tikzcd}[cramped]
	K & H \\
	{K\cap L} & {H\cap L}
	\arrow["{\not\in\cT}"', dotted, from=1-1, to=1-2]
	\arrow["{\cS}", from=2-1, to=1-1]
	\arrow["\cT", dashed, from=2-1, to=2-2]
	\arrow["{\cS}"', from=2-2, to=1-2]
\end{tikzcd} \]
    
    Since $K \cap (H \cap L)=K\cap L$, the above square contradicts the second requirement of Definition \ref{def:compatiblePair}. Thus $(\cS, \cT)$ cannot be a compatible pair of transfer systems.
\end{proof}

\begin{theorem}
    \label{thm:maximalLinear}
    Let $G$ be an abelian group and let $\cU$ be a $G$-universe, then $\Hull(\cT_{\cD(\cU)})=\cT_{\cL(\overline\cU)}$.
\end{theorem}
\begin{proof}
    From definition \ref{def:universeClosure}, we have $\cT_{\cD(\cU)}=\cT_{\cD(\overline\cU)}$. Thus it suffices to show that $\Hull(\cT_{\cD(\cU)})=\cT_{\cL(\cU)}$ for a closed universe $\cU$. This is proven by a sequence of refinements $\Hull(\cT_{\cD(\cU)})\leq \cT_{\cL(\cU)} \leq \Hull(\cT_{\cD(\cU)})$. The second refinement follows from the fact that $(\cT_{\cL(\cU)},\cT_{\cD(\cU)})$ is a compatible pair of transfer systems, see \cite[Proposition 6.16]{blumberg2021} or \cite[Theorem 5.1]{chan2025realizingcompatiblepairstransfer}. For the first inequality, assume that $K \to_{\Hull(\cD(\cU))} H$. By Theorem \ref{thm:internalLinearisHull}, we have an refinement $\coind^H_K \res^G_K \cT_{\cD(\cU)} \leq \res^G_H \cT_{\cD(\cU)}$. Since $\cD$ commutes with restriction and takes induction to coinduction, this implies that $\cT_{\cD(\ind^H_K \res^G_K\cU)} \leq \cT_{\cD(\res^G_H \cU)}$. The $G$-universe $\cU$ is closed, hence by Lemma \ref{lem:resClosed}, the $H$-universe $\res^G_H\cU$ is closed. By Lemma \ref{lem:closedProperties}, $\ind^H_K \res^G_K \cU$ embeds into $\res^G_H \cU$ by $H$-equivariant linear maps, meaning $K \to_{\cL(\cU)} H$. 
\end{proof}

In light of the results of \cite{chan2025realizingcompatiblepairstransfer}, we have the following corollary.

\begin{corollary}
    \label{cor:realizingCompatiblePairs}
    If $G$ is an abelian group and $\cU$ is a $G$-universe, then for any compatible pair of transfer systems $(\cT,\cT_{\cD(\cU)})$, there exists an $N_\infty$-operad pair $(\cP,\cO)$ such that $\cT=\cT_\cP$ and $\cO \simeq \cD(\cU)$.
\end{corollary}
\begin{proof}
    By \cite[Proposition 7.16]{chan2025realizingcompatiblepairstransfer}, it is sufficient to show that the maximal compatible pair of transfer systems $(\Hull(\cT_{\cD(\cU)}),\cT_{\cD(\cU)})$ can be realised as an operad pair. By Theorem \ref{thm:maximalLinear}, this compatible pair is realised by the Steiner operad pair $(\cL(\overline{\cU}),\cK(\overline\cU))$.
\end{proof}

\section{Linear isometries operads and saturated transfer systems}
\label{section:saturated}
We have seen that for all universes $\cU$, the transfer system $\cT_{\cL(\cU)}$ is saturated. By direct computation with groups of small order, Rubin \cite{rubin2021detecting} proved that saturated transfer systems need not come from linear isometries operads in general. Following this, the question of determining for which finite groups saturated transfer systems always come from linear isometries operads has seen much attention, particularly in the abelian case, see \cite{rubin2021detecting,hafeez2022saturatedlinearisometrictransfer,macbrough2023equivariant,xiaequivariant}. 

\begin{definition}
    \label{def:saturatedgroup}
    A finite group $G$ is \textbf{saturated} if, for each saturated $G$ transfer system $\cT$, there exists a $G$ universe $\cU$ such that $\cT = \cT_{\cL(\cU)}$.
\end{definition}

For a prime $p$, Rubin proved that all cyclic groups of order $p^n$ are saturated in \cite{rubin2021detecting}. Later, MacBrough showed that cyclic groups of non-prime-power order are saturated as long as their order is coprime to $6$, and that all abelian groups of rank at least $3$ are not saturated in \cite{macbrough2023equivariant}. Thus, the only remaining abelian groups $G$ where it is not known in general whether $G$ is saturated are those of rank $2$. 

The non-abelian case has not seen quite so much success. Apart from a few non-abelian groups of small order, which have been shown not to be saturated in \cite[\S 5]{rubin2021detecting}, Xia proved that non-abelian $p$-groups $G$ are not saturated if they have subgroups $H$ and $K$ satisfying particular criteria in \cite[Theorem 3.4]{xiaequivariant}. No saturated non-abelian groups have yet been found in the current literature. Given the evidence, we propose the following conjecture.

\begin{conjecture}
    \label{conj:saturatedAbelian}
    If $G$ is a saturated finite group, then $G$ is abelian.
\end{conjecture}

We will outline a general strategy for proving non-abelian groups are not saturated, building on \cite{xiaequivariant}. In particular, we will apply the results of this section to all non-abelian symmetric, alternating and Hamiltonian groups, and describe situations where our current arguments fail. We begin with the following observation.

\begin{proposition}
    \label{prop:completeUniverse}
    Let $\cU$ be a universe (either real or complex) such that $\{e\} \to_{\cL(\cU)} G$, then $\cU$ is a complete universe. 
\end{proposition}
\begin{proof}
    The real case is \cite[Corollary 5.5]{rubin2021detecting}, and the complex case can be proven identically. 
\end{proof}

Let $G$ be a non-abelian finite group, and $G'$ a subgroup of $G$. By Proposition \ref{prop:changeofgroupUniverses}, if $\cU$ is a $G$-universe, then $\res^G_{G'} \cT_{\cL(\cU)} =\cT_{\cL(\res^G_{G'}\cU)}$. In particular, if a saturated transfer system $\cT$ can be realised as a linear isometries operad $\cL(\cU)$, then $\res^G_{G'} \cT$ is realised by $\cL(\res^G_{G'} \cU)$. Our strategy is to find a small subgroup $G'$ of $G$, and a $G$ transfer system $\cT$ such that $\res^G_{G'} \cT$ cannot be realised as a linear isometries operad. It will then follow that $\cT$ also does not arise from any linear isometries operad. Since the non-abelian group $G$ is assumed to be finite in this article, it has finitely many non-abelian subgroups $G'$. Choosing a non-abelian subgroup $G'$ which is minimal with respect to subgroup inclusion (among non-abelian subgroups of $G$) gives a group in which every proper subgroup is abelian. Such groups are known as minimal non-abelian groups. These will be the subgroups that we are restricting to. 
  
In \cite[Corollary 3.5]{xiaequivariant}, Xia proves that for all odd primes $p$, any finite minimal non-abelian group of order $p^\alpha$ is not saturated. Using the following results of Miller-Moreno, we will extend Xia's argument to all minimal non-abelian finite groups.

\begin{theorem}[\cite{miller1903non}]
    \label{thm:miller}
    Let $G$ be a minimal non-abelian finite group. Then the following statements hold. 
    \begin{enumerate}
        \item The order of $G$ has at-most $2$ prime factors $p$ and $q$.
        \item $G$ has a normal subgroup of index $p$.
        \item If $|G|=p^\alpha q^\beta$ for $\alpha,\beta \geq 1$, then: 
        \begin{itemize}
            \item There is one subgroup of order $q^\beta$.
            \item There are $q^\beta$ subgroups of order $p^\alpha$. 
        \end{itemize} 
    \end{enumerate} 
\end{theorem}

In particular, we need only consider groups whose order can be divided by at most two primes. As in \cite{macbrough2023equivariant} and \cite{xiaequivariant}, we will work with complex representations. The benefit of complex representations here is that an irreducible complex $G$-representation $V$ is one-dimensional if and only if the commutator $G^{(1)}$ acts trivially on $V$. We define $\fI_{>1}$ to be the complex $G$-universe containing precisely all of the irreducible representations of degree greater than one. Thus $\fI_{>1}$ contains all representations on which $G^{(1)}$ acts nontrivially. This property can be leveraged to prove the following lemmas, which are analogous to \cite[Lemmas 3.2 \& 3.3]{xiaequivariant}.

\begin{lemma}[{\cite[Lemma 3.3]{xiaequivariant}}]
    \label{lem:completeOnK}
    Let $G$ be a finite group and suppose that there exists a subgroup $K \leq G$ such that $K$ does not contain the commutator $G^{(1)}$ of $G$. If $\fU$ is a complex $G$-universe such that $\fI_{>1}$ embeds into $\fU$, then $\res^G_K \fU $ is a complete $K$-universe. 
\end{lemma}

\begin{remark}
    While $G$ is explicitly assumed to be a $p$-group for an odd prime $p$ throughout \cite{xiaequivariant}, this assumption is not used in the proof of \cite[Lemma 3.3]{xiaequivariant}. The same is not true for \cite[Lemma 3.2]{xiaequivariant}, and we shall explain how the argument generalises below. We warn the reader that \cite{xiaequivariant} uses subsets of $\Irr_\C(G)$ in place of complex $G$-universes. However, these form isomorphic lattices, under which $\fI_{>1}$ corresponds to $\Irr_\C^{>1}(G)$.
\end{remark}

\begin{lemma}
    \label{lem:completeOnH}
    Let $G$ be a finite group of order $p^\alpha q^\beta$, where $p$ and $q$ are distinct primes, $\alpha \geq 1$ and $\beta \geq 0$. Suppose that there exists a normal subgroup $H \unlhd G$ such that $H$ is abelian and $[G: H]=p$. Then, if $\fU$ is a complex $G$-universe such that $\res^G_H\fU$ is a complete $H$-universe, it follows that $\fI_{>1}$ embeds into $\fU$. 
\end{lemma}
\begin{proof}
   The prime power case of this lemma is treated in \cite[Lemma 3.2]{xiaequivariant}. The present lemma can be proven mostly analogously, although more subtlety is required for a few of the steps. 
   
   Let $|G^{(1)}|=p^\gamma q^\delta$, then we obtain the following two formulae (Compare \cite[Lemma 3.2]{xiaequivariant}).

    \begin{equation}
        \label{formula:completeOnHequality}
        \sum_{V \in \Irr_\C^{>1}(G)}(\dim V)^2 = |G|-|G_{\mathrm{ab}}|=p^\alpha q^\beta - p^{\alpha - \gamma }q^{\beta-\delta}.
    \end{equation}
    
    \begin{equation}
        \label{formula:completeOnHinequality}
        \sum_{V \in \Irr^{>1}_\C (G)} \dim V \geq p^{\alpha-1}q^\beta - p^{\alpha -1 -\gamma }q^{\beta - \delta}.
    \end{equation}

    As in the proof of \cite[Lemma 3.2]{xiaequivariant}, the right side of Inequality \ref{formula:completeOnHinequality} is the number of irreducible $H$-representations $W$ such that $W$ has a non-trivial action by $G^{(1)}$. Each such $W$ must be one-dimensional since $H$ is abelian. The left side represents the fact that $W$ must embed into $\res^G_H V$, for some irreducible $G$-representation $V$ of dimension greater than one. 
    
    Recall that inductions of representations of finite groups are ambidextrous in the sense that $\ind^G_H \dashv \res^G_H \dashv \ind^G_H$. Thus, for all irreducible $H$-representations $W$ as in the above paragraph, the irreducible subrepresentations of $\ind^G_H W$ are precisely the irreducible $G$-representations $V$, such that $W$ embeds into $\res^G_H V$. The dimension of $V$ is greater than one and must also divide the order of $G$. In particular, we have that $p \mid\dim V$ or $q \mid \dim V$. Note that $\dim (\ind^G_H W)=[G:H]=p$, so if $p \mid \dim V$, then $\dim V = p$ and $V = \ind^G_H W$. Otherwise $q \mid \dim V$ for all irreducible subrepresentations $V$ of $\ind^G_H W$. However, this implies that $q \mid \dim (\ind^G_H W)=p$, contradicting the fact that $p$ and $q$ are distinct primes. Thus, the dimension of $V$ must be $p$. 

    We can conclude analogously to the proof of \cite[Lemma 3.2]{xiaequivariant}, that both sides of Inequality \ref{formula:completeOnHinequality} divide Equation \ref{formula:completeOnHequality} by a factor of $p$. Hence \ref{formula:completeOnHinequality} must, in fact, be an equality. Thus, in order for $\res^G_H \fU$ to contain each irreducible $H$-representation $W$, the $G$-universe $\fU$ must contain each irreducible $G$-representation $V$ of dimension greater than one.
\end{proof}

Putting Lemmas \ref{lem:completeOnK}  and \ref{lem:completeOnH} together, we can prove the analogue of \cite[Theorem 3.4]{xiaequivariant}.

\begin{theorem}
    \label{thm:pqgroupUnsaturated}
    Let $G$ be a finite group of order $p^\alpha q^\beta$ for $\alpha \geq 1$ and $\beta \geq 0$, and suppose that: 
    \begin{enumerate}
        \item There is an abelian normal subgroup $H\unlhd G$ of index $p$.
        \item There is a subgroup $K \leq G$ such that $G^{(1)}\nleq K$.
    \end{enumerate}
    Then, for all (real) $G$-universes $\cU$, if $\{e\} \to_{\cL(\cU)} H$, then $\{e\} \to_{\cL(\cU)} K$.
\end{theorem}
\begin{proof}
    Let $\cU$ be a real $G$-universe such that $\{e\} \to_{\cL(\cU)} H$. If $\fU$ is the complexification of $\cU$, then $\cT_{\cL(\fU)}=\cT_{\cL(\cU)}$. Since $\{e\} \to_{\cL(\fU)} H$ it follows, by Proposition \ref{prop:completeUniverse}, that $\res^G_H\fU$ must be a complete $H$-universe. Thus, we can apply Lemma \ref{lem:completeOnH} to deduce that  $\fI_{>1}$ embeds into $\fU$. Lemma \ref{lem:completeOnK} then tells us that $\res^G_K \fU$ is a complete $K$-universe. We can conclude that $\{e\}\to_{\cL(\cU)}K$.
\end{proof}

In particular, if $G,H,K$ are groups as in the above theorem, and $\cT$ is a saturated $G$ transfer system where $\{e\} \to_\cT H$ holds but $\{e\} \to_\cT K$ does not, then $\cT$ cannot be realised as a linear isometries operad.  
We shall apply Theorem \ref{thm:pqgroupUnsaturated} to minimal non-abelian groups. First, consider the case that $|G|$ has exactly two prime factors.

\begin{lemma}
    \label{lem:pqalmostAbelian}
    Let $G$ be a minimal non-abelian group, and suppose that $|G|=p^\alpha q^\beta $ for two distinct primes $p$ and $q$ and two integers $\alpha,\beta > 0$. Then there exist subgroups $H,K \leq G$ such that the following statements (or the equivalent statements with $p$ and $q$ swapped) hold:
    \begin{itemize}
        \item $H$ is normal of index $p$.
        \item $[G:K]=q^\beta$ and $G^{(1)}\nleq K$.
    \end{itemize}
    In particular, any group which is isomorphic to $H$ cannot contain $K$ as a subgroup.  
\end{lemma}
\begin{proof}
   By Theorem \ref{thm:miller}, $G$ has a normal subgroup $H$ of index $p$. Now take $K$ to be any one of the subgroups of order $p^\alpha$. Again, by Theorem \ref{thm:miller}, we know that there are precisely $q^\beta$ such subgroups. By the Sylow theorems, $K$ is conjugate to all other subgroups of order $p^\alpha$, thus it cannot be normal. Therefore, $K$ does not contain the commutator subgroup $G^{(1)}$.  
\end{proof}

We shall use the following construction as a counterexample to certain non-abelian groups being saturated.

\begin{definition}
    For a group $G$ and a subgroup $H \leq G$, define the saturated transfer system $\cT_H^G$ by the rule that $L \to_{\cT_H^G} K$ if either $L=K$ or there exists $g \in G$ such that $L \leq K \leq H^g$.
\end{definition}

\begin{remark}
    Using the induction functor of Section \ref{subsec:changeGroup}, we may alternatively write $\cT_H^G$ as $\ind^G_H \cT_H$ where $\cT_H$ denotes the maximal $H$-transfer system.
\end{remark}

\begin{theorem}
    \label{thm:nonsaturated1}
    Let $G$ be a non-abelian group. Suppose that there is a minimal non-abelian subgroup $G'\leq G$ of order $|G'| = p^\alpha q^\beta$ for $\alpha,\beta > 0$. Then $G$ is not saturated.
\end{theorem}
\begin{proof}
    Since $G'$ is minimal non-abelian of order $|G'|=p^\alpha q^\beta$ for $\alpha,\beta > 0$, we can find subgroups $H, K \leq G'$ as in Lemma \ref{lem:pqalmostAbelian}. We will show that the saturated transfer system $\cT_H^G$ cannot be realised as a linear isometries operad. 

    Suppose that there exists a $G$-universe $\cU$ such that $\cT_H^G=\cT_{\cL(\cU)}$. We shall proceed by seeking a contradiction. By construction, we have a transfer $\{e\} \to_{\cL(\res^G_{G'} \cU)} H$, which implies that $\{e\} \to_{\cL(\res^G_{G'}\cU)} K$ by Theorem \ref{thm:pqgroupUnsaturated}. Using the results of Section \ref{subsec:changeGroup}, we have $\cT_{\cL(\res^G_{G'} \cU)}=\res^G_{G'}\cT_{\cL(\cU)}=\res^G_{G'} \cT^G_H$. Thus, there is a transfer $\{e\} \to_{\cT^G_H} K$. By definition, this implies that $K \leq H^g$ for some $g \in G$. However, Lemma \ref{lem:pqalmostAbelian} tells us that $H$ and $K$ have mutually indivisible orders, giving the desired contradiction. It follows that the saturated $G$ transfer system $\cT^G_H$ cannot be realised by a linear isometries operad.
\end{proof}

As an immediate consequence of Theorem \ref{thm:nonsaturated1}, we can classify all of the saturated symmetric and alternating groups.

\begin{corollary}
\label{cor:symmetricNotsaturated}
The symmetric group $\Sigma_n$ is saturated if and only if $n\leq 2$, and the alternating group $A_m$ is saturated if and only if $m \leq 3$.  
\end{corollary}
\begin{proof}
    Note that the groups $\Sigma_1$, $A_1$ and $A_2$ are all trivial, and thus trivially saturated. The groups $\Sigma_2$ and $A_3$ are cyclic of prime order and are therefore saturated. 
    
    Now consider $\Sigma_3$. We have $|\Sigma_3|=6$, $\Sigma_3$ is not abelian and each proper subgroup of $\Sigma_3$ is cyclic. Thus, by Theorem \ref{thm:nonsaturated1}, any group which contains $\Sigma_3$ as a subgroup is not saturated. For $n \geq 3$, the subgroup $\langle(123),(12) \rangle\leq \Sigma_n$ is isomorphic to $\Sigma_3$, hence $\Sigma_n$ is not saturated. 
    Similarly, $|A_4|=12$ and $A_4$ contains only abelian proper subgroups, hence any group which contains $A_4$ as a subgroup is not saturated. For $m \geq 4$, the subgroup $\langle (123),(12)(34)\rangle \leq A_m$ is isomorphic to $A_4$, so $A_m$ is not saturated.
\end{proof}

We will now examine the prime power case for minimal non-abelian groups. Odd primes were covered by Xia in \cite[Corollary 3.5]{xiaequivariant}. We will expand Xia's results to the case $p=2$ and apply these results to show that Hamiltonian groups are never saturated.

\begin{theorem}
    Let $p$ be a prime and let $P$ be a minimal non-abelian group of order $p^\alpha$, then $P$ is not saturated.
\end{theorem}
\begin{proof}
      For odd primes $p$, the present lemma is proven in \cite[Corollary 3.5]{xiaequivariant}, using the prime power variant of Lemma \ref{lem:pqalmostAbelian}. However, the assumption that $p$ is odd is only necessary to deduce that $P$ contains a non-cyclic subgroup. By \cite[\S 5 Theorem 4.10]{gorenstein2007finite}, a $p$-group with no non-cyclic abelian subgroups is either cyclic or generalised Quaternion. If $P$ is not a generalised Quaternion group, then the proof of the present lemma is identical to the proof of \cite[Corollary 3.5]{xiaequivariant}. Thus, it remains to address the generalised Quaternion groups $Q_{4n}$ with $n \geq 2$. Since $Q_{4n}$ contains $Q_8$ as a subgroup, the only minimal non-abelian generalised Quaternion group is $Q_8$ itself. Thus, the lemma can be deduced from Rubin's computations in \cite[Example 5.10]{rubin2021detecting}.
\end{proof}

It might seem at first like we have all of the pieces in place to prove Conjecture \ref{conj:saturatedAbelian}. After all, every non-abelian group contains a minimal non-abelian group, and we have shown that no minimal non-abelian group is saturated. However, for groups $G' \leq G$, it is not true in general that a saturated $G'$-transfer system is in the image of the restriction map $\res^G_{G'}$, even for transfer systems of the form $\cT^{G'}_H$ for a prime-index normal subgroup $H$ of $G'$. In the case that $G'$ was a minimal non-abelian group of order $|G'|=p^\alpha q^\beta$, this problem could be avoided. This was because we could construct $H$ and $K$ in Theorem \ref{thm:pqgroupUnsaturated} such that the orders of $H$ and $K$ were mutually indivisible. As such, a saturated $G$ transfer system containing the arrow $\{e\} \to H$ but not $\{e\} \to K$ could easily be constructed. This argument immediately fails in the prime-power case, where there are no two subgroups with mutually indivisible orders. Despite this, we can still find applications in the case that the subgroup inclusion we are restricting along admits a retract.

\begin{lemma}
    \label{lem:splitInclusion}
    Let $G$ be a saturated group and let $H$ be a subgroup such that the inclusion $\iota:H \to G$ admits a retract $\pi$, then $H$ is saturated. 
\end{lemma}
\begin{proof}
    Let $\cT$ be a saturated $H$-transfer system. By functoriality, we have that $\iota^{-1}_L \circ \pi^{-1}_L = \mathrm{id}_{\trans(H)}$. The map $\iota^{-1}_L$ is just $\res^G_H$, and $\pi^{-1}_L\cT$ is the inflation of $\cT$ along the epimorphism $\pi$. By Proposition \ref{prop:changeDiscSaturated}, $\cT$ is equal to the restriction of the saturated $G$ transfer system $\pi^{-1}_L\cT$. Since $G$ is saturated, this transfer system can be realised as a linear isometries operad $\cL(\cU)$, thus $\cT$ can be realised by $\cL(\res^G_H\cU)$.
\end{proof}

\begin{corollary}
\label{cor:hamiltonianNotSaturated}
    If $G$ is a finite Hamiltonian group (i.e. $G$ is non-abelian, but every subgroup of $G$ is normal), then $G$ is not saturated.
\end{corollary}
\begin{proof}
    It is known that every finite Hamiltonian group $G$ splits as a direct product of $Q_8$ and an abelian group $A$. Thus, the inclusion $Q_8 \to G$ admits a retract, given by the usual projection map. We know that $Q_8$ is not saturated. Hence, by Lemma \ref{lem:splitInclusion}, $G$ is not saturated. 
\end{proof}

\bibliographystyle{alpha}
\bibliography{reference_items.bib}

\end{document}